\documentclass[twoside,11pt]{article}

\pdfoutput=1

\usepackage[english]{babel}
\usepackage{amsmath}
\usepackage{amsthm}
\usepackage{amsfonts}
\usepackage{amssymb}
\usepackage{graphicx}
\usepackage{mathrsfs}

\allowdisplaybreaks[3]

 \theoremstyle{plain}
\newtheorem{theorem}{Theorem}[section]
\theoremstyle{remark}
\newtheorem{remark}[theorem]{\bf Remark}

\numberwithin{equation}{section}

\newcommand{\be}{\begin{equation}}
\newcommand{\ee}{\end{equation}}

\begin{document}

\title{{\LARGE Low Mach number limit of strong solutions to the compressible primitive equations with gravity}
 \footnotetext{\small E-mail addresses: mupc024@163.com
 }
 }

\author{{Pengcheng Mu}\\[2mm]
\small\it School of Mathematical Sciences, Anhui University, Hefei 230601, P.R.China
}

\date{}

\maketitle

\begin{quote}
\small \textbf{Abstract}: In this paper, we explore the low Mach number singular limit of the local-in-time strong solutions to the compressible primitive equations with gravity for general adiabatic coefficient. First we construct the uniform estimate for the solutions to the non-dimensional compressible primitive equations with general ill-prepared initial data. Due to the effects of gravity and the anisotropy of the system, the operator with large coefficient in this model is not explicitly skew-symmetric. Thus, obtaining the uniform estimate requires novel techniques. After that, we investigate rigorously the low Mach number limit of the compressible primitive equations with both well-prepared and ill-prepared initial data. The limiting system is shown to be the incompressible primitive equations with inhomogeneous density that depends on the vertical variable.


\indent \textbf{AMS (2010) Subject Classification}: 35Q35; 76U05; 76M45
\end{quote}

\section{Introduction}\label{sec1}

We consider the following 3-D isentropic rotating compressible primitive equations describing the motion of flows in geophysics:
\begin{equation}\label{CPE}
\begin{cases}
\partial_{t}\rho+\mathrm{div}_{h}\,(\rho u)+\partial_{z}(\rho w)=0,\\
\partial_{t}(\rho u)+\mathrm{div}_{h}\,(\rho u\otimes u)+\partial_{z}(\rho u w)+f\rho u^{\bot}+\nabla_{h}P(\rho)\\
~~~~~~~~~~=\mu\triangle_{h}u+\mu\partial_{z}^{2}u+\nu\nabla_{h}\mathrm{div}_{h}u,\\
\partial_{z}P(\rho)=-g\rho.
\end{cases}
\end{equation}
In \eqref{CPE}, the unknown $\rho$ denotes the density of the fluid, $u=(u_{1},u_{2})^{\top}$ the horizontal velocity field, and $w$ the vertical velocity. The pressure $P$ is given by $P(\rho)=\rho^{\gamma}$ where $\gamma>1$ is the adiabatic coefficient. The positive constants $f$, $g$ denote the rotating frequency and the gravitational acceleration, respectively, and $\mu$, $\nu$ are the viscosity coefficients. The independent variables in \eqref{CPE} are $t\in\mathbb{R}^{+}$, $x=(x_{1},x_{2})\in \mathbb{T}_{L\times L}^{2}$ and $z\in(0,L_{z})$, where $\mathbb{T}_{L\times L}^{2}=[0,2\pi L]\times [0,2\pi L]$ is a 2-D torus, and $L_{z}$ is the height of the domain.  $\nabla_{h}=(\partial_{x_{1}},\partial_{x_{2}})^{\top}$, $\mathrm{div}_{h}=\nabla_{h}\cdot$ denote the horizontal gradient and divergence, respectively, and $\triangle_{h}=\partial_{x_{1}}^{2}+\partial_{x_{2}}^{2}$ is the horizontal Laplace.

It is well-known that the planetary horizontal scales are significantly large than the vertical one in geophysical flows in the large scale motion. Using this property the authors in \cite{B7121} formally derived the compressible primitive equations from the compressible Navier-Stokes equations. Different from the 3-D isentropic compressible Navier-Stokes equations, the equation of the vertical velocity (vertical momentum) in \eqref{CPE} is replaced by the hydrostatic balance equation, i.e., \eqref{CPE}$_{3}$. This equation describes the balance between the vertical variation of the pressure with the force created by gravity, and is used as a fundamental equation in geophysical models \cite{B6221,B7122,B7123}.

These years, there arise some mathematical researches related to the well-possdness and asymptotics of the compressible primitive equations. First, in \cite{B921,B922} the authors constructed the global weak solutions to the simplified 2-D compressible primitive equations. Then, in \cite{B7127} Tang-Gao showed the stability of weak solutions to the 2-D compressible primitive equations by using the entropy estimates. Later, in 2018, Jiu-Li-Wang \cite{B923} proved that the weak solution constructed in \cite{B922} is unique. In 2021, Wang \cite{B925} showed the existence and uniqueness of the local strong solution to the 2-D compressible primitive equations with vacuum in the initial data. For the 3-D compressible primitive equations, Ersoy-Ngom-Sy \cite{B7121} studied the stability of weak solutions to \eqref{CPE} with density-dependent viscosity after it was formally derived from the 3-D Navier-Stokes equations. After that, by using the methods in \cite{B7125,B7126}, the authors in \cite{B6221,B7101} independently constructed the global weak solutions of compressible primitive equations \eqref{CPE} with $f=0$ (without rotation), $g=0$ (without gravity), and density-dependent viscosity. Later, in 2021, Liu-Titi \cite{B7128} established the local-in-time well-posedness of strong solutions to \eqref{CPE} with $f=0$ in the following two cases: with gravity but without vacuum, and with vacuum but without gravity. In 2023, Tang-Ne\v{c}asov\'{a} rigorously derived the inviscid compressible primitive equations from the Euler system by using the relative entropy inequality \cite{B9151}.

There are also some results concerning the low Mach number limit of the compressible primitive equations in the past several years. In 2020, Gao-Ne\v{c}asov\'{a}-Tang \cite{B941} proved the weak-strong uniqueness of \eqref{CPE} and studied the low Mach number limit of weak solutions to \eqref{CPE} with well-prepared initial data. Liu-Titi \cite{B71121} investigated the low Mach number limit of local strong solutions to \eqref{CPE} constructed in \cite{B7128}. By elaborately estimating the difference between \eqref{CPE} and the limiting system, the authors in \cite{B71121} showed that the strong solutions of \eqref{CPE} with well-prepared initial data converge to that of the incompressible primitive equations as the Mach number tends to zero. The convergence rate and a class of global large solutions to \eqref{CPE} for small Mach number are also obtained. Later, in 2023, the low Mach number limit for the local strong solutions to \eqref{CPE} with ill-prepared initial data was also studied by the same authors in \cite{B944}. We mention here that, the works in \cite{B941,B71121,B944} above all neglected the effects of gravity. Indeed, considering gravity in \eqref{CPE} will cause challenging difficulties \cite{B71121}. The only work so far concerning the singular limit problem of the compressible primitive equations with gravity, to our best knowledge,  was done in \cite{Dai-Mu}. In that paper, the authors considered the low Mach number limit of the compressible primitive equations \eqref{CPE} with $g>0$ for a special case $\gamma=2$, and showed that the local strong solutions of the system with well-prepared initial data, as well as their time derivatives, converge to those of the inhomogeneous incompressible primitive equations. However, for general $\gamma>1$, the low Mach number limit of \eqref{CPE} with gravity is still not solved. In this work, we shall address this problem.

More precisely, the aim of the present paper is to explore rigorously the low Mach number singular limit of local-in-time strong solutions to \eqref{CPE} {\it with gravity for all $\gamma>1$}. Indeed, different from \cite{B941,B71121,B944}, when gravity is taken into account, the density of the fluid is no longer $z$-independent. This fact, combined with the anisotropy of the equations, makes the operator with large coefficient (i.e., the reciprocal of the Mach number) not explicitly skew-symmetric. Thus, it's troublesome to construct the uniform estimates for the solutions. Although the equations \eqref{CPE} can be reformulated into a form in which the large operator is skew-symmetric after being multiplied by some weighted function in the case $\gamma=2$ in \cite{Dai-Mu}, this property no longer holds in this paper. Indeed, when $\gamma\neq2$ an extra penalized term (i.e, $\delta^{-1}(\gamma-2)\theta w$ in \eqref{eq2}$_{1}$) appears in the density equation after it is reformulated. This term can not be annihilated by performing usual estimations. Thus, when deriving the uniform bounds for the solutions, we adopt a new refined estimating strategy which is different from that in \cite{Dai-Mu} and shown to be more subtle.

Before studying the low Mach number limit of the system \eqref{CPE}, we write it in the following non-dimensional form:
\begin{equation}\label{e1}
\begin{cases}
\partial_{t}\rho+\mathrm{div}_{h}\,(\rho u)+\partial_{z}(\rho w)=0,\\
\rho(\partial_{t}u+u\cdot\nabla_{h} u+w\cdot\partial_{z}u+\frac{1}{Ro}u^{\bot})+\frac{1}{\gamma\delta^{2}}\nabla_{h}\rho^{\gamma}=\mathbb{S}u,\\
\partial_{z}\rho^{\gamma-1}=-\theta.
\end{cases}
\end{equation}
The derivations from the original system \eqref{CPE} to the non-dimensional form \eqref{e1} are presented in appendix. In \eqref{e1}, the parameters $Ro$, $\delta$ are the Rossby and Mach numbers, respectively, and $\theta$ is a positive constant. The non-dimensional density $\rho$ and viscosity $\mathbb{S}u$ are given by
\begin{equation}\label{eqrho}
\rho=(1-\theta z+\delta r)^{\frac{1}{\gamma-1}},
\end{equation}
and
\begin{equation}\label{eqS}
\mathbb{S}u=\mu\triangle_{h}u+\mu\partial_{z}^{2}u+\nu\nabla_{h}\mathrm{div}_{h}\,u,
\end{equation}
where $r$ is the unknown fluctuation of the density in horizontal directions such that $\partial_{z}r=0$, and $\mu$, $\nu$ are non-dimensional viscosity coefficients which are assumed to be fixed constants in this paper. The non-dimensional equations \eqref{e1} are set in the domain
\begin{equation*}
x=(x_{1},x_{2})\in\mathbb{T}^{2}:=[0,2\pi]\times[0,2\pi],~~0<z<\kappa:=\frac{L_{z}}{L},
\end{equation*}
where $\kappa$ is known as the aspect ratio parameter which is usually small in large scale motion in geophysical science. In this paper, we are not taking $\kappa\to0$. However, to avoid the vacuum, we assume that
\begin{equation*}
\inf\limits_{z}\,(1-\theta z)=1-\theta\kappa\geq\frac{1}{2},~~\mathrm{thus}~~\theta\kappa\leq\frac{1}{2}.
\end{equation*}
By introducing \eqref{eqrho} into \eqref{e1}, dividing \eqref{e1}$_{2}$ by $\rho$ and denoting by
\begin{equation}\label{eqXi}
\xi:=1-\theta z+\delta r,
\end{equation}
we get
\begin{equation}\label{eq2}
\begin{cases}
\partial_{t}r+\frac{\gamma-1}{\delta}\mathrm{div}_{h}\,(\xi u)+\frac{\gamma-1}{\delta}\partial_{z}(\xi w)-(\gamma-2)u\cdot\nabla_{h}r+\frac{(\gamma-2)\theta}{\delta}w=0,\\
\partial_{t}u+u\cdot\nabla_{h} u+w\cdot\partial_{z}u+\frac{1}{Ro}u^{\bot}+\frac{1}{(\gamma-1)\delta}\nabla_{h}r=\frac{1}{\rho}\mathbb{S}u,\\
\partial_{z}r=0.
\end{cases}
\end{equation}
It is clear that systems \eqref{eq2} and \eqref{e1} are equivalent for smooth solutions that keep away from vacuum. We impose \eqref{eq2} with initial data
\begin{equation}\label{in}
(r,u)|_{t=0}=(r_{0},u_{0})
\end{equation}
and the stress-free boundary condition
\begin{equation}\label{bd}
(\partial_{z}u,w)|_{z=0,\kappa}=0.
\end{equation}
In this paper, we shall first construct the uniform estimates for strong solutions to the problem \eqref{eq2}-\eqref{in}-\eqref{bd} in Section \ref{sec2}. After that, we then study the low Mach number singular limit of \eqref{eq2} with well-prepared initial data, and derive explicitly the convergence rate as the Mach number tends to zero in Section \ref{sec3}. Besides, in Section \ref{sec4}, the low Mach number limit of \eqref{eq2} with ill-prepared initial data shall also be considered. Finally, we present the derivations from the original compressible primitive equations \eqref{CPE} to the non-dimensional system \eqref{e1} in appendix.

Now we state the main results of this paper. We denote by
\begin{equation}\label{eq-nota}
\Omega:=\mathbb{T}^{2}\times(0,\kappa),~~~\xi^{L}:=1-\theta z,~~~\mathrm{and}~~~\rho^{L}:=(\xi^{L})^{\frac{1}{\gamma-1}}=(1-\theta z)^{\frac{1}{\gamma-1}}.
\end{equation}

\begin{theorem}\label{th1}$\mathrm{(Uniform\,estimates).}$
Assume that
\begin{equation}\label{eq8253}
\theta\kappa\leq
\begin{cases}
\min\Big\{\frac{1}{2},\frac{\sqrt{2}(\gamma-1)}{8|\gamma-2|}\Big\},~~\mathrm{if}~~\gamma>1,~\gamma\neq2;\\
\frac{1}{2},~~~~\mathrm{if}~~\gamma=2.
\end{cases}
\end{equation}
Assume further that the initial data $(r_{0},u_{0})$ satisfy
\begin{equation}\label{eq974}
\|u_{0}\|_{H^{2}(\Omega)}+\|\partial_{h}^{3}u_{0}\|_{L^{2}(\Omega)}+\|r_{0}\|_{H^{3}(\mathbb{T}^{2})}\leq B_{in},~~\partial_{z}u_{0}|_{z=0,\kappa}=0,
\end{equation}
where $B_{in}$ is the upper bound of the initial data and is independent of $\delta$. Then there exists $\delta_{0}>0$ such that for any $0<\delta\leq \delta_{0}$, the problem \eqref{eq2}-\eqref{in}-\eqref{bd} admits a unique solution $(r,u)$ such that
\begin{equation*}
\begin{split}
&r\in L^{\infty}(0,T;H^{3}(\mathbb{T}^{2})),~~~u\in L^{\infty}(0,T;H^{2}(\Omega)),~~~\partial_{h}^{3}u\in L^{\infty}(0,T;L^{2}(\Omega)),\\
&(\nabla_{h}u,\partial_{z}u)\in L^{2}(0,T;H^{2}(\Omega)),~~~(\nabla_{h}\partial_{h}^{3}u,\partial_{z}\partial_{h}^{3}u)\in L^{2}(0,T;L^{2}(\Omega)),
\end{split}
\end{equation*}
and it holds
\begin{equation}\label{est}
\begin{split}
&\sup\limits_{t\in[0,T]}\big(\|r(t)\|^{2}_{H^{3}(\mathbb{T}^{2})}+\|u(t)\|^{2}_{H^{2}(\Omega)}+\|\partial_{h}^{3}u(t)\|^{2}_{L^{2}(\Omega)}\big)+\int_{0}^{T}{\|(\nabla_{h}u,\partial_{z}u)(t)\|^{2}_{H^{2}(\Omega)}}dt\\
&+\int_{0}^{T}{\|(\nabla_{h}\partial_{h}^{3}u,\partial_{z}\partial_{h}^{3}u)(t)\|^{2}_{L^{2}(\Omega)}}dt\leq C\big(\|r_{0}\|^{2}_{H^{3}(\mathbb{T}^{2})}+\|u_{0}\|^{2}_{H^{2}(\Omega)}+\|\partial_{h}^{3}u_{0}\|^{2}_{L^{2}(\Omega)}\big),
\end{split}
\end{equation}
where the positive time $T$ and the constant $C$ are independent of $\delta$.
\end{theorem}

For any $\delta>0$, the local well-posedness of the strong solution to \eqref{e1} (thus \eqref{eq2}) on some time interval $[0,T_{\delta}]$ can be proved by adapting the discussions in \cite{B7128} when the initial density is away from vacuum. So, to prove Theorem \ref{th1}, it suffices to show that $T_{\delta}\geq T$ for some positive time $T$ which is independent of $\delta$, and the uniform estimate \eqref{est} holds.
However, this is not an easy task since the equations \eqref{eq2} are not explicitly skew-symmetric in the sense that the penalized term $\delta^{-1}(\gamma-2)\theta w$ in \eqref{eq2}$_{1}$ can not be annihilated by performing the usual estimations. Thus, to get rid of this difficulty, we shall first derive the uniform estimates for the horizontal derivatives of $(\rho,u)$ from the equations \eqref{e1}, and then show that these estimates are equivalent to those for $(r,u)$ in \eqref{est}
in a short time when $\delta$ is small enough. We mention here that, the hydrostatic balance equation \eqref{e1}$_{3}$ plays an important role in our proof, although it do bring some troubles in this problem.

\begin{remark}
The estimates for up to the second derivatives of $(r,u)$ are not enough to deal with the integrals, for example, \eqref{eq9231} and \eqref{eq9232} in Section \ref{sec2}. Thus, in \eqref{est}, we need extra estimates for the third horizontal derivatives. Indeed, these estimates are utilized to bound the $L^{\infty}(0,T;H^{2}(\Omega))$ norm of the vertical velocity $w$. Noting that, from \eqref{eq921}, the expression of $w$ involves an integral over the vertical variable, thus the estimates for the vertical derivatives of $u$ are less required than that for the horizontal ones in the estimations of $w$. This is why the estimate for $\|\partial_{z}^{3}u\|_{L^{\infty}(0,T;L^{2}(\Omega))}$ is not involved in \eqref{est}.
\end{remark}

\begin{remark}
From \eqref{est} and the estimates for $w$ in Section \ref{sec2} we can also get, under the conditions in \eqref{eq8253}, that
\begin{equation}\label{eq-est-w}
\|w\|_{L^{\infty}(0,T;H^{2}(\Omega))}\leq C.
\end{equation}
\end{remark}

\begin{remark}
From the proof in Section \ref{sec2} we shall see that, besides the Mach number, the estimates in \eqref{est} are also independent of the Rossby number $Ro$ due to the skew symmetry of the rotation term. Thus, the asymptotics of \eqref{eq2} as the Rossby number tends to zero can also be studied.
\end{remark}

Next, as the Mach number $\delta$ tends to zero, formally, the velocity $u$ in \eqref{eq2} shall converge to the solution of the following 3-D inhomogeneous incompressible primitive equations:
\begin{equation}\label{eqL}
\begin{cases}
\partial_{t}u^{L}+u^{L}\cdot\nabla_{h}u^{L}+w^{L}\cdot\partial_{z}u^{L}+\frac{1}{Ro}(u^{L})^{\bot}+\nabla_{h}\Pi=\frac{1}{\rho^{L}}\mathbb{S}u^{L},\\
\mathrm{div}_{h}(\rho^{L}u^{L})+\partial_{z}(\rho^{L}w^{L})=0,
\end{cases}
\end{equation}
where $\Pi=\Pi(x,t)$ is the Lagrange multiplier such that $\partial_{z}\Pi=0$ and $\int_{\mathbb{T}^{2}}{\Pi(x,t)}dx=0$, and $\rho^{L}$ and $\mathbb{S}$ are defined by \eqref{eq-nota} and \eqref{eqS}, respectively. In the following, we will justify this limit for \eqref{eq2} for both well-prepared and ill-prepared initial data. Hereafter, we use
\begin{equation*}
\mathbf{P}:=\mathrm{id}-\nabla_{h}\triangle_{h}^{-1}\mathrm{div}_{h}
\end{equation*}
to denote the Leray operator in $\mathbb{T}^{2}$, and for any $F:\Omega\mapsto \mathbb{R}$, we define
\begin{equation*}
\overline{F}:=\frac{1}{\kappa}\int^{\kappa}_{0}{F}dz,~~\mathrm{and}~~\widetilde{F}:=F-\overline{F}.
\end{equation*}

\begin{theorem}\label{th2}$\mathrm{(Low\,Mach\,number\,limit\,with\,well\mbox{-}prepared\,data).}$
Let $\{(r,u)\}_{\delta}$ be a sequence of local strong solutions to the problem \eqref{eq2}-\eqref{in}-\eqref{bd} given by Theorem \ref{th1}. If
\begin{equation}\label{eq983}
\|r_{0}\|_{H^{1}(\mathbb{T}^{2})}+\|u_{0}-u_{0}^{L}\|_{H^{1}(\Omega)}\leq C\delta,~~\mathrm{div}_{h}\,(\overline{\rho^{L}u^{L}_{0}})=0,~~\mathrm{and}~~u_{0}^{L}\in H^{3}(\Omega),
\end{equation}
then we have
\begin{equation}\label{eqWell}
\sup\limits_{t\in[0,T]}\big(\|r(t)\|^{2}_{H^{1}(\mathbb{T}^{2})}+\|(u-u^{L})(t)\|^{2}_{H^{1}(\Omega)}\big)+\int_{0}^{T}{\|(\nabla_{h},\partial_{z})(u-u^{L})(t)\|^{2}_{H^{1}(\Omega)}}dt\leq C\delta^{2},
\end{equation}
where the constant $C$ above is independent of $\delta$, and $u^{L}$ is the strong solution of the system \eqref{eqL}
supplemented with the following initial and boundary conditions:
\begin{equation}\label{inL}
u^{L}|_{t=0}=u^{L}_{0},~~~(\partial_{z}u^{L},w^{L})|_{z=0,\kappa}=0.
\end{equation}

\end{theorem}

\begin{theorem}\label{th3}$\mathrm{(Low\,Mach\,number\,limit\,with\,ill\mbox{-}prepared\,data).}$
Let $\{(r,u)\}_{\delta}$ be a sequence of local strong solutions to the problem \eqref{eq2}-\eqref{in}-\eqref{bd} given by Theorem \ref{th1} and $\rho$ be defined by \eqref{eqrho}.
Assume that
\begin{equation}\label{eq9251}
u_{0}\to u^{L}_{0}~~\mathrm{weakly~in}~~H^{2}(\Omega),~~\mathrm{as}~~\delta\to0,
\end{equation}
then as $\delta\to0$,
\begin{equation*}
\rho\to\rho^{L}~~\mathrm{strongly~in}~~L^{\infty}(0,T;H^{3}(\Omega)),
\end{equation*}
\begin{equation*}
(u,w)\to(u^{L},w^{L})~~\mathrm{wakely*~in}~~L^{\infty}(0,T;H^{2}(\Omega)),
\end{equation*}
and
\begin{equation*}\
(\partial_{z}u,\mathrm{curl}_{h}u)\to(\partial_{z}u^{L},\mathrm{curl}_{h}u^{L})~~\mathrm{strongly~in}~~L^{\infty}(0,T;L^{2}(\Omega)),
\end{equation*}
where $(\rho^{L},u^{L},w^{L})$ is the strong solution of the equations \eqref{eqL} supplemented with the following initial and boundary conditions:
\begin{equation}\label{inLL}
u^{L}|_{t=0}=\mathbf{P}u^{L}_{0},~~~(\partial_{z}u^{L},w^{L})|_{z=0,\kappa}=0.
\end{equation}

\end{theorem}

The difference between well-prepared and ill-prepared initial data is that the well-prepared data have excluded the acoustic waves in the system, while the ill-prepared data allow the propagation of the fast-oscillating acoustic waves in the solutions. As a result, in Theorem \ref{th2}, we can show the strong convergence of the solutions and derive explicitly the convergence rate by performing energy estimates for the difference between the systems \eqref{eq2} and \eqref{eqL}, with more refined treatments on terms involved with pressure. However, when the data are ill-prepared, the solutions to \eqref{eq2} are expected to converge weakly as $\delta\to0$ in $\Omega$ which is periodic in horizontal directions. We compare the results and proof of Theorem \ref{th3} with those in \cite{B944} where the authors explored the low Mach number limit of the compressible primitive equations without gravity for ill-prepared initial data. Indeed, when gravity is neglected, the density is $z$-independent thus the compressible primitive system endowed with stress-free boundary condition can be extended as a purely periodic problem with respect to the vertical variable. As a result, the spectrum method based on Fourier transform can be applied to deal with the resonances of the fast oscillating waves in the periodic domain, or the dispersions in whole space. Nevertheless, in this paper, as the density depends on $z$, \eqref{e1} can not be recast as a periodic system in the vertical direction. Thus, the spectrum method in \cite{B944} is inapplicable here. Actually, the proof of Theorem \ref{th3} involves passing to the limit $\delta\to0$ in the weak formulation (i.e, the integral formula) of equations \eqref{eq2}, after obtaining necessary compactness from \eqref{est}. The strong convergences of $\mathrm{curl}_{h}u$ and $\partial_{z}u$ are key in the identification of the limits of convection terms.

\begin{remark}
Noting that, different from the incompressible primitive equations with constant density in \cite{B71121,B944}, the density $\rho^{L}$ in \eqref{eqL} depends on the vertical variable $z$. However, since $2^{-1/(\gamma-1)}\leq \rho^{L}\leq 1$, the local well-posedness of the strong solution to \eqref{eqL} with initial data in $H^{2}(\Omega)$ can also be established by repeating the discussions in  \cite{B7112,B7113}. Moreover, by adapting the proof in \cite{B7114,B71121}, one can further obtain the global regularity for \eqref{eqL} under some constraints on the viscosity coefficients.
\end{remark}


\begin{remark}
In similar fashion to \cite{Dai-Mu}, when the data are well-prepared, the strong convergence of the time derivatives to $(r,u)$ can also be proved if more conditions are imposed on $r_{0}$, $u_{0}$ and $u^{L}_{0}$.
\end{remark}

We end this section by appointing some notations that shall be used. Throughout this paper, we use $\partial_{h}^{n}$ to denote the $n$-th horizontal spatial derivatives, and $C$ the positive constants that are independent of $\delta$ and may take different values in each appearance. The notation $A\approx B$ means that $A\leq CB$ and meanwhile $B\leq CA$ for some `irrelevant' constant $C$. Furthermore, we shall use the following inequalities
\begin{equation}\label{eq971}
\|F\|_{L^{3}(\Omega)}\leq C\|F\|_{L^{2}(\Omega)}^{1/2}\|F\|_{H^{1}(\Omega)}^{1/2},~~\|F\|_{L^{\infty}(\Omega)}\leq C\|F\|_{L^{2}(\Omega)}^{1/4}\|F\|_{H^{2}(\Omega)}^{3/4}.
\end{equation}
\begin{equation}\label{eq8271}
\|\int^{z}_{0}{F}(x,\tau)d\tau\|_{L^{2}(\Omega)}\leq \frac{1}{\sqrt{2}}\kappa\|F\|_{L^{2}(\Omega)},~~\|\overline{F}\|_{L^{p}(\mathbb{T}^{2})}+\|\widetilde{F}\|_{L^{p}(\Omega)}\leq C\|F\|_{L^{p}(\Omega)},~~\forall\,p\in[1,\infty].
\end{equation}
In particular, for the special case $p=2$ we have
\begin{equation}\label{eq9114}
\kappa\|\overline{F}\|^{2}_{L^{2}(\mathbb{T}^{2})}+\|\widetilde{F}\|^{2}_{L^{2}(\Omega)}= \|F\|^{2}_{L^{2}(\Omega)}.
\end{equation}

\section{Uniform estimates}\label{sec2}
In this section we construct the uniform estimates for the strong solutions to \eqref{eq2}. For the initial data $(r_{0},u_{0})$ given by Theorem \ref{th1}, we let
\begin{equation}\label{eq8262}
\delta_{*}:=\frac{1}{8C_{*}B_{in}},
\end{equation}
where $C_{*}$ is the continuity modulus of $H^{2}(\mathbb{T}^{2})\hookrightarrow L^{\infty}(\mathbb{T}^{2})$. Then when $\delta\leq \delta_{*}$ we have
\begin{equation*}
\rho_{0}:=\Big(1-\theta z+\delta r_{0}\Big)^{\frac{1}{\gamma-1}}\geq\Big(1-\theta\kappa-\delta\|r_{0}\|_{L^{\infty}(\mathbb{T}^{2})}\Big)^{\frac{1}{\gamma-1}}\geq \Big(\frac{1}{2}-\delta_{*} C_{*}B_{in}\Big)^{\frac{1}{\gamma-1}}\geq\Big(\frac{3}{8}\Big)^{\frac{1}{\gamma-1}},
\end{equation*}
i.e., the initial density is away from vacuum. Thus, for any $\delta\leq \delta_{*}$, the unique local strong solution $(\rho,u)$ to \eqref{e1} with boundary condition \eqref{bd} and initial data $(\rho_{0},u_{0})$ on some time interval $[0,T_{\delta}]$ can be constructed by adapting the steps in \cite{B7128}, and $\rho$ is away from vacuum on that interval. Thus, the problem \eqref{eq2}-\eqref{in}-\eqref{bd} is also well posed on $[0,T_{\delta}]$. In the following, we shall show that $T_{\delta}$ admit a uniform lower bound $T>0$ which is independent of $\delta$ as $\delta$ is small enough, and the estimate \eqref{est} holds.

First, we introduce the notations:
\begin{equation}\label{eq8251}
\mathcal{A}_{*}=\mathcal{A}_{*}(t):=\big(\|u\|_{H^{2}(\Omega)}+\|\partial_{h}^{3}u\|_{L^{2}(\Omega)}+\|r\|_{H^{3}(\mathbb{T}^{2})}\big)(t),
\end{equation}
\begin{equation*}
\begin{split}
\mathcal{A}=\mathcal{A}(t)&:=\Big(\int_{\Omega}{\rho|u|^{2}+\rho|\partial_{h}u|^{2}+\rho|\partial_{h}^{2}u|^{2}+|\partial_{z}u|^{2}+|\partial_{z}^{2}u|^{2}+|\partial_{z}\partial_{h}u|^{2}+\rho|\partial_{h}^{3}u|^{2}}dxdz\\
&+\int_{\Omega}{\frac{2}{\gamma(\gamma-1)\delta^{2}}H(\rho)+\frac{1}{\delta^{2}}\rho^{\gamma-2}|\partial_{h}\rho|^{2}+\frac{1}{\delta^{2}}\rho^{\gamma-2}|\partial_{h}^{2}\rho|^{2}+\frac{1}{\delta^{2}}\rho^{\gamma-2}|\partial_{h}^{3}\rho|^{2}}dxdz\Big)^{1/2}(t),
\end{split}
\end{equation*}
and
\begin{equation*}
\begin{split}
\mathcal{B}=\mathcal{B}(t):=&\Big(-\sum\limits_{n=0}^{3}\int_{\Omega}{\mathbb{S}(\partial_{h}^{n}u)\cdot\partial_{h}^{n}u}dxdz-\int_{\Omega}{\mathbb{S}(\partial_{z}u)\cdot(\partial_{z}u-\partial_{z}^{3}u)}dxdz\\
&-\int_{\Omega}{\mathbb{S}(\partial_{h}\partial_{z}u)\cdot\partial_{h}\partial_{z}u}dxdz\Big)^{1/2}(t),
\end{split}
\end{equation*}
where
\begin{equation}\label{eq82220}
H(\rho)=\rho^{\gamma}-\gamma\xi^{L}\rho+(\gamma-1)\times(\xi^{L})^{\frac{\gamma}{\gamma-1}}
\end{equation}
with $\xi^{L}$ defined by \eqref{eq-nota}. By an application of integration by parts, it is clear that
\begin{equation*}
\mathcal{B}\approx \|\nabla_{h}u\|_{H^{2}(\Omega)}+\|\partial_{z}u\|_{H^{2}(\Omega)}+\|\nabla_{h}\partial_{h}^{3}u\|_{L^{2}(\Omega)}+\|\partial_{z}\partial_{h}^{3}u\|_{L^{2}(\Omega)}.
\end{equation*}
Now we show that $\mathcal{A}$ and $\mathcal{A}_{*}$ are equivalent when $t$ and $\delta$ are small enough.

Indeed, by direct calculations we have
\begin{equation}\label{eq8229}
\begin{split}
&\partial_{h}\rho=\frac{\delta}{\gamma-1}\rho^{2-\gamma}\partial_{h}r,~~~~\partial_{h}^{2}\rho=\frac{\delta}{\gamma-1}\rho^{2-\gamma}\partial_{h}^{2}r+\frac{(2-\gamma)\delta^{2}}{(\gamma-1)^{2}}\rho^{3-2\gamma}\cdot|\partial_{h}r|^{2},~~~~\mathrm{and}\\
&\partial_{h}^{3}\rho=\frac{\delta}{\gamma-1}\rho^{2-\gamma}\cdot\partial_{h}^{3}r+\frac{3(2-\gamma)\delta^{2}}{(\gamma-1)^{2}}\rho^{3-2\gamma}\cdot\partial_{h}r\cdot\partial_{h}^{2}r+\frac{(2-\gamma)(3-2\gamma)\delta^{3}}{(\gamma-1)^{3}}\rho^{4-3\gamma}\cdot|\partial_{h}r|^{3}.
\end{split}
\end{equation}
Thus
\begin{equation}\label{eq973}
\begin{split}
&\int_{\Omega}{\frac{1}{\delta^{2}}\rho^{\gamma-2}|\partial_{h}^{3}\rho|^{2}}dxdz=C\int_{\Omega}{\rho^{2-\gamma}}|\partial_{h}^{3}r|^{2}dxdz\\
&~~~~+\int_{\Omega}{C\delta^{2}\rho^{4-3\gamma}|\partial_{h}r|^{2}\cdot|\partial_{h}^{2}r|^{2}+C\delta^{4}\rho^{6-5\gamma}|\partial_{h}r|^{6}+C\delta\rho^{3-2\gamma}\partial_{h}^{3}r\cdot\partial_{h}r\cdot\partial_{h}^{2}r}dxdz\\
&~~~~+\int_{\Omega}{C\delta^{2}\rho^{4-3\gamma}\partial_{h}^{3}r\cdot|\partial_{h}r|^{3}+C\delta^{3}\rho^{5-4\gamma}\partial_{h}^{2}r\cdot|\partial_{h}r|^{4}}dxdz
\end{split}
\end{equation}
for some different constants $C$ that only depend on $\gamma$. As a result, if
\begin{equation}\label{eq8228}
\Big(\frac{1}{4}\Big)^{\frac{1}{\gamma-1}}\leq \rho\leq \Big(\frac{5}{4}\Big)^{\frac{1}{\gamma-1}},
\end{equation}
then we can verify that
\begin{equation*}
\begin{split}
&\int_{\Omega}{C\delta^{2}\rho^{4-3\gamma}|\partial_{h}r|^{2}\cdot|\partial_{h}^{2}r|^{2}+C\delta^{4}\rho^{6-5\gamma}|\partial_{h}r|^{6}+C\delta\rho^{3-2\gamma}\partial_{h}^{3}r\cdot\partial_{h}r\cdot\partial_{h}^{2}r}dxdz\\
&~~~+\int_{\Omega}{C\delta^{2}\rho^{4-3\gamma}\partial_{h}^{3}r\cdot|\partial_{h}r|^{3}+C\delta^{3}\rho^{5-4\gamma}\partial_{h}^{2}r\cdot|\partial_{h}r|^{4}}dxdz\\
&\leq C\big(\delta\|r\|_{H^{3}(\mathbb{T}^{2})}+\delta^{2}||r\|^{2}_{H^{3}(\mathbb{T}^{2})}+\delta^{3}\|r\|_{H^{3}(\mathbb{T}^{2})}^{3}+\delta^{4}\|r\|_{H^{3}(\mathbb{T}^{2})}^{4}\big)\|r\|_{H^{3}(\mathbb{T}^{2})}^{2}.\\
\end{split}
\end{equation*}
Introducing the above inequality into \eqref{eq973} we get
\begin{equation}\label{eq8221}
\begin{split}
&\int_{\Omega}{\frac{1}{\delta^{2}}\rho^{\gamma-2}|\partial_{h}^{3}\rho|^{2}}dxdz\lesseqgtr C\|\partial_{h}^{3}r\|_{L^{2}(\mathbb{T}^{2})}^{2}\pm C\Big(\sum\limits_{n=1}^{4}\big(\delta\|r\|_{H^{3}(\mathbb{T}^{2})}\big)^{n}\Big)\|r\|_{H^{3}(\mathbb{T}^{2})}^{2}.\\
\end{split}
\end{equation}
In similar fashion, we can check
\begin{equation}\label{eq8222}
\int_{\Omega}{\frac{2}{\gamma(\gamma-1)\delta^{2}}H(\rho)}dxdz \approx \frac{1}{\delta^{2}}\|\rho-(\xi^{L})^{\frac{1}{\gamma-1}}\|_{L^{2}(\Omega)}^{2}\approx \|r\|_{L^{2}(\mathbb{T}^{2})}^{2},
\end{equation}
\begin{equation}\label{eq8224}
\begin{split}
\int_{\Omega}{\frac{1}{\delta^{2}}\rho^{\gamma-2}|\partial_{h}\rho|^{2}}dxdz \approx \|\partial_{h}r\|_{L^{2}(\mathbb{T}^{2})}^{2},
\end{split}
\end{equation}
and
\begin{equation}\label{eq8223}
\begin{split}
\int_{\Omega}{\frac{1}{\delta^{2}}\rho^{\gamma-2}|\partial_{h}^{2}\rho|^{2}}dxdz\lesseqgtr& C\|\partial_{h}^{2}r\|_{L^{2}(\mathbb{T}^{2})}^{2}\pm C\big(\delta\|r\|_{H^{3}(\mathbb{T}^{2})}+\delta^{2}\|r\|_{H^{3}(\mathbb{T}^{2})}^{2}\big)\|r\|_{H^{3}(\mathbb{T}^{2})}^{2}.
\end{split}
\end{equation}
Combining \eqref{eq8221}, \eqref{eq8222}, \eqref{eq8224} and \eqref{eq8223} we conclude that, if \eqref{eq8228} holds, then there exist some positive constants $C_{1}$, $C_{2}$, $C_{3}$ and $C_{4}$, which only depend on $\gamma$, such that
\begin{equation}\label{eq8225}
\begin{split}
&C_{1}\|r\|_{H^{3}(\mathbb{T}^{2})}^{2}-C_{2}\Big(\sum\limits_{n=1}^{4}\big(\delta\|r\|_{H^{3}(\mathbb{T}^{2})}\big)^{n}\Big)\|r\|_{H^{3}(\mathbb{T}^{2})}^{2}\\
&\leq\int_{\Omega}{\frac{2}{\gamma(\gamma-1)\delta^{2}}H(\rho)+\frac{1}{\delta^{2}}\rho^{\gamma-2}|\partial_{h}\rho|^{2}+\frac{1}{\delta^{2}}\rho^{\gamma-2}|\partial_{h}^{2}\rho|^{2}+\frac{1}{\delta^{2}}\rho^{\gamma-2}|\partial_{h}^{3}\rho|^{2}}dxdz\\
&\leq C_{3}\|r\|_{H^{3}(\mathbb{T}^{2})}^{2}+ C_{4}\Big(\sum\limits_{n=1}^{4}\big(\delta\|r\|_{H^{3}(\mathbb{T}^{2})}\big)^{n}\Big)\|r\|_{H^{3}(\mathbb{T}^{2})}^{2}.
\end{split}
\end{equation}
Due to the above estimation, if \eqref{eq8228} and
\begin{equation}\label{eq8252}
C_{2}\Big(\sum\limits_{n=1}^{4}\big(\delta\|r\|_{H^{3}(\mathbb{T}^{2})}\big)^{n}\Big)\leq \frac{1}{2}C_{1}
\end{equation}
both hold true we can further get
\begin{equation*}
\begin{split}
\int_{\Omega}{\frac{2}{\gamma(\gamma-1)\delta^{2}}H(\rho)+\frac{1}{\delta^{2}}\rho^{\gamma-2}|\partial_{h}\rho|^{2}+\frac{1}{\delta^{2}}\rho^{\gamma-2}|\partial_{h}^{2}\rho|^{2}+\frac{1}{\delta^{2}}\rho^{\gamma-2}|\partial_{h}^{3}\rho|^{2}}dxdz\approx \|r\|_{H^{3}(\mathbb{T}^{2})}^{2},
\end{split}
\end{equation*}
which immediately implies
\begin{equation}\label{eq82210}
\frac{1}{C_{5}}\mathcal{A}_{*}\leq \mathcal{A}\leq C_{6}\mathcal{A}_{*}
\end{equation}
for some constants $C_{5}$ and $C_{6}$ that are independent of $\delta$. We assume that $C_{5}C_{6}>1$, otherwise we can take $C_{6}=C_{6}+\frac{2}{C_{5}}$.

Now we define
\begin{equation}\label{eq8261}
T^{*}_{\delta}:=\sup\limits_{t}\Big\{t\leq T_{\delta}:\sup\limits_{s\in[0,t]}\mathcal{A}_{*}^{2}(s)+\frac{C_{5}^{2}}{2}\int^{t}_{0}{\mathcal{B}^{2}(s)}ds\leq 4C_{5}^{2}C_{6}^{2}\mathcal{A}^{2}_{*}(0)\Big\},
\end{equation}
where $C_{5}$ and $C_{6}$ are the constants in \eqref{eq82210}. Since $4C_{5}^{2}C_{6}^{2}>1$,
by the continuity of the solution, we have $T^{*}_{\delta}>0$. Then, we take a positive constant $b_{0}$ such that
\begin{equation*}
C_{2}\Big(\sum\limits_{n=1}^{4}b_{0}^{n}\Big)\leq\frac{1}{2}C_{1},
\end{equation*}
and let
\begin{equation}\label{eq9253}
\delta_{**}:=\min\Big\{\delta_{*},\frac{1}{8C_{*}C_{5}C_{6}B_{in}},\frac{b_{0}}{2C_{5}C_{6}B_{in}}\Big\},
\end{equation}
where $C_{1}$ and $C_{2}$ are from \eqref{eq8225}, and $\delta_{*}$, $C_{*}$ and $B_{in}$ are from \eqref{eq8262}. Then when $t\in(0,T^{*}_{\delta}]$ and $\delta\leq \delta_{**}$ we have
\begin{equation*}
\delta\|r\|_{L^{\infty}(\mathbb{T}^{2})}\leq \delta_{**}C_{*}\mathcal{A}_{*}\leq \frac{1}{8C_{*}C_{5}C_{6}B_{in}}\times C_{*}\times2C_{5}C_{6}B_{in}\leq \frac{1}{4},
\end{equation*}
and
\begin{equation*}
C_{2}\Big(\sum\limits_{n=1}^{4}\big(\delta\|r\|_{H^{3}(\mathbb{T}^{2})}\big)^{n}\Big)\leq C_{2}\Big(\sum\limits_{n=1}^{4}\big(\delta_{**}\mathcal{A}_{*}\big)^{n}\Big)\leq C_{2}\Big(\sum\limits_{n=1}^{4}b_{0}^{n}\Big)\leq \frac{1}{2}C_{1},
\end{equation*}
where we have used the fact $\mathcal{A}_{*}(0)\leq B_{in}$ from \eqref{eq974}. Thus we get \eqref{eq8228} and \eqref{eq8252} by using the condition in \eqref{eq8253}. As a result, from the above discussions we conclude that \eqref{eq82210} also holds.

In the following we restrict $t\in(0,T^{*}_{\delta}]$ and $\delta\leq \delta_{**}$, and use \eqref{eq8229}, \eqref{eq8228}, \eqref{eq82210} and the definition of $\mathcal{A}$ to derive that
\begin{equation}\label{eq82225}
\begin{split}
&\|\partial_{h}^{2}\rho\|_{L^{4}(\Omega)}\leq C\delta\big(\mathcal{A}+\mathcal{A}^{2}\big),~~~~~\|\rho\|_{H^{2}(\Omega)}\leq C\big(1+\delta\mathcal{A}\big),\\
&\|\partial_{z}\partial_{h}^{3}\rho\|_{L^{2}(\Omega)}\leq C\big(\delta\|\partial_{h}^{3}r\|_{L^{2}(\Omega)}+\delta^{2}\|r\|^{2}_{H^{3}(\Omega)}+\delta^{3}\|r\|^{3}_{H^{3}(\Omega)}\big)\leq C\delta\big(\mathcal{A}+\mathcal{A}^{2}+\mathcal{A}^{3}\big),\\
&\|\partial_{h}^{2}\rho^{\gamma}\|_{L^{2}(\Omega)}=\|\partial_{h}^{2}\rho\cdot\rho^{\gamma-1}+2\partial_{h}\rho\cdot\partial_{h}\rho^{\gamma-1}+\rho\cdot\partial_{h}^{2}\rho^{\gamma-1}\|_{L^{2}(\Omega)}\leq C\delta\big(\mathcal{A}+\mathcal{A}^{2}\big),\\
&\|\partial_{h}\rho\|_{L^{\infty}(\Omega)}\leq C\delta\|r\|_{H^{3}(\Omega)}\leq C\delta\mathcal{A},~~~~~\mathrm{and}~~~~~\|\partial_{z}\rho\|_{L^{\infty}(\Omega)}\leq C.
\end{split}
\end{equation}
These estimates shall be frequently used later.

\textbf{Estimates~for~$w$.} To derive the estimations for $w$, we first average \eqref{eq2}$_{1}$ over $(0,\kappa)$ to obtain
\begin{equation}\label{eq8142}
\partial_{t}r+\frac{\gamma-1}{\delta}\mathrm{div}_{h}\,(\overline{\xi u}) -(\gamma-2)\overline{u}\cdot\nabla_{h}r+\frac{(\gamma-2)\theta}{\delta}\overline{w}=0.
\end{equation}
Subtracting \eqref{eq8142} from \eqref{eq2}$_{1}$ we get
\begin{equation*}
(\gamma-1)\mathrm{div}_{h}\,(\widetilde{\xi u})+(\gamma-1)\partial_{z}(\xi w)-\delta(\gamma-2)\widetilde{u}\cdot\nabla_{h}r+(\gamma-2)\theta\widetilde{w}=0,
\end{equation*}
which implies that
\begin{equation}\label{eq921}
w(z)=-\frac{1}{\xi}\int^{z}_{0}{\Big(\mathrm{div}_{h}\,(\widetilde{\xi u})-\frac{\delta(\gamma-2)}{\gamma-1}\widetilde{u}\cdot\nabla_{h}r\Big)(\tau)}d\tau-\frac{(\gamma-2)\theta}{(\gamma-1)\xi}\int^{z}_{0}{\widetilde{w}(\tau)}d\tau.
\end{equation}
Then, using \eqref{eq921}, \eqref{eq8253}, \eqref{eq8271}, \eqref{eq9114} and \eqref{eq8228} we have
\begin{equation*}
\begin{split}
&\|w\|_{L^{2}(\Omega)}\leq C\|\xi^{-1}\|_{L^{\infty}(\Omega)}\|(\mathrm{div}_{h}\,(\xi u),\delta u\cdot\nabla_{h}r)\|_{L^{2}(\Omega)}+\frac{|\gamma-2|}{\gamma-1}\theta\|\xi^{-1}\|_{L^{\infty}(\Omega)}\|\int^{z}_{0}{\widetilde{w}}d\tau\|_{L^{2}(\Omega)}\\
&\leq C\big(\|\xi\|_{L^{\infty}(\Omega)}\|\mathrm{div}_{h}\,u\|_{L^{2}(\Omega)}+\delta\|\nabla_{h}r\|_{L^{4}(\Omega)}\|u\|_{L^{4}(\Omega)}\big)+\frac{|\gamma-2|}{\sqrt{2}(\gamma-1)}\theta\kappa\|\xi^{-1}\|_{L^{\infty}(\Omega)}\|w\|_{L^{2}(\Omega)}\\
&\leq C\big(1+\delta\|r\|_{H^{2}(
\Omega)}\big)\|u\|_{H^{2}(
\Omega)}+\frac{4|\gamma-2|}{\sqrt{2}(\gamma-1)}\theta\kappa\|w\|_{L^{2}(\Omega)}\\
&\leq C\big(1+\mathcal{A}\big)\mathcal{A}+\frac{1}{2}\|w\|_{L^{2}(\Omega)},
\end{split}
\end{equation*}
which leads to
\begin{equation*}
\begin{split}
\|w\|_{L^{2}(\Omega)}\leq C\big(\mathcal{A}+\mathcal{A}^{2}\big).
\end{split}
\end{equation*}
In similar fashion, we can verify
\begin{equation*}
\begin{split}
&\|\partial_{h}w\|_{L^{2}(\Omega)}\leq C\delta\|\partial_{h}r\|_{L^{\infty}(\Omega)}\|(\mathrm{div}_{h}\,(\xi u),\delta u\cdot\nabla_{h}r,w)\|_{L^{2}(\Omega)}\\
&~~~~+C\|(\partial_{h}\mathrm{div}_{h}\,(\xi u),\delta\partial_{h}(u\cdot\nabla_{h}r))\|_{L^{2}(\Omega)}+\frac{|\gamma-2|}{\gamma-1}\theta\|\xi^{-1}\|_{L^{\infty}(\Omega)}\|\int^{z}_{0}{\partial_{h}\widetilde{w}}d\tau\|_{L^{2}(\Omega)}\\
&\leq C\delta\|r\|_{H^{3}(\Omega)}\Big(\big(1+\delta\|r\|_{H^{2}(\Omega)}\big)\|u\|_{H^{2}(\Omega)}+\|w\|_{L^{2}(\Omega)}\Big)+C\big(1+\delta\|r\|_{H^{2}(\Omega)}\big)\|u\|_{H^{2}(\Omega)}\\
&~~~+\frac{4|\gamma-2|}{\sqrt{2}(\gamma-1)}\theta\kappa\|\partial_{h}w\|_{L^{2}(\Omega)}\\
&\leq C\mathcal{A}^{2}\big(1+\mathcal{A}\big)+C\big(1+\mathcal{A}\big)\mathcal{A}+\frac{1}{2}\|\partial_{h}w\|_{L^{2}(\Omega)},
\end{split}
\end{equation*}
and
\begin{equation*}
\begin{split}
\|\partial_{z}w\|_{L^{2}(\Omega)}&\leq C\big(\|\xi^{-2}\|_{L^{\infty}(\Omega)}+\|\xi^{-1}\|_{L^{\infty}(\Omega)}\big)\|(\mathrm{div}_{h}\,(\xi u),\delta u\cdot\nabla_{h}r,w)\|_{L^{2}(\Omega)}\\
&\leq C\big(1+\delta\|r\|_{H^{2}(
\Omega)}\big)\|u\|_{H^{2}(
\Omega)}+C\|w\|_{L^{2}(\Omega)}\leq C\big(\mathcal{A}+\mathcal{A}^{2}\big),
\end{split}
\end{equation*}
thus
\begin{equation}\label{eq82230}
\|w\|_{H^{1}(\Omega)}\leq C\big(\mathcal{A}+\mathcal{A}^{2}+\mathcal{A}^{3}\big) \leq C\big(\mathcal{A}+\mathcal{A}^{3}\big).
\end{equation}
Finally, by using similar methods we can also check that
\begin{equation}\label{eq8232}
\|w\|_{H^{2}(\Omega)}\leq C\big(\mathcal{A}+\mathcal{A}^{4}\big).
\end{equation}

\begin{remark}
From the above discussions we can see that, different from the case $\gamma=2$ in \cite{Dai-Mu}, the vertical velocity $w$ in \eqref{eq921} can not be explicitly expressed by the horizontal velocities and the density. Thus, to derive the estimates for $w$,
we require the condition $\theta\kappa\leq \sqrt{2}(\gamma-1)/(8|\gamma-2|)$ in \eqref{eq8253}.
\end{remark}

\textbf{Estimates~for~$(r,u)$.} Now we begin to construct the uniform estimates for $(r,u)$ and their derivatives. First, by taking the $L^{2}$ inner product of \eqref{e1}$_{2}$ with $u$, integrating by parts and taking use of \eqref{e1}$_{1,3}$ we obtain with
\begin{equation}\label{eq8122}
\frac{d}{dt}\int_{\Omega}{\frac{1}{2}\rho|u|^{2}+\frac{1}{\gamma(\gamma-1)\delta^{2}}\rho^{\gamma}}dxdz+\frac{\theta}{(\gamma-1)\delta^{2}}\int_{\Omega}{\rho w}dxdz-\int_{\Omega}{\mathbb{S}u\cdot u}dxdz=0.
\end{equation}
The integral involved with the vertical velocity $w$ above is troublesome to us because the coefficient of which is singular as $\delta\to0$, and the sign of $\rho w$ is actually unknown. To deal with this integral, the key new observation here is the following identity:
\begin{equation}\label{eq82212}
\partial_{t}(\xi^{L}\rho)+\mathrm{div}_{h}\,(\xi^{L}\rho u)+\partial_{z}(\xi^{L}\rho w)+\theta\rho w=0.
\end{equation}
This identity is obtained by multiplying \eqref{e1}$_{1}$ by $\xi^{L}$ an using the fact $\partial_{z}\xi^{L}=-\theta$. Then, integrating \eqref{eq82212} over $\Omega$ we acquire
\begin{equation}\label{eq8311}
\frac{d}{dt}\int_{\Omega}{\xi^{L}\rho}dxdz+\theta\int_{\Omega}{\rho w}dxdz=0.
\end{equation}
The above equality can help us to annihilate the third term in \eqref{eq8122}, as well as to obtain the desired estimate for $H(\rho)$. Indeed, from \eqref{eq8311} we can verify that
\begin{equation}\label{eq8133}
\begin{split}
&\frac{d}{dt}\int_{\Omega}{\frac{1}{\gamma(\gamma-1)\delta^{2}}\rho^{\gamma}}dxdz+\frac{\theta}{(\gamma-1)\delta^{2}}\int_{\Omega}{\rho w}dxdz\\
&=\frac{d}{dt}\int_{\Omega}{\frac{1}{\gamma(\gamma-1)\delta^{2}}H(\rho)}dxdz+\frac{1}{(\gamma-1)\delta^{2}}\frac{d}{dt}\int_{\Omega}{\xi^{L}\rho}dxdz+\frac{\theta}{(\gamma-1)\delta^{2}}\int_{\Omega}{\rho w}dxdz\\
&=\frac{d}{dt}\int_{\Omega}{\frac{1}{\gamma(\gamma-1)\delta^{2}}H(\rho)}dxdz,
\end{split}
\end{equation}
where $H(\rho)$ is defined by \eqref{eq82220}. Thus, plugging \eqref{eq8133} into \eqref{eq8122} we finally get
\begin{equation}\label{eqs2}
\begin{split}
&\frac{1}{2}\frac{d}{dt}\int_{\Omega}{\bigg(\rho|u|^{2}+\frac{2}{\gamma(\gamma-1)\delta^{2}}H(\rho)\bigg)}dxdz-\int_{\Omega}{\mathbb{S}u\cdot u}dxdz=0.
\end{split}
\end{equation}

\textbf{Estimates~for~$(\partial_{h}r,\partial_{h}u)$.} Applying $\partial_{h}$ to \eqref{e1}$_{1,2}$ we get
\begin{equation}\label{eq8134}
\begin{cases}
\partial_{t}\partial_{h}\rho+\mathrm{div}_{h}\,(\rho\partial_{h}u)+\mathrm{div}_{h}\,(\partial_{h}\rho\cdot u)+\partial_{z}\partial_{h}(\rho w)=0,\\
\rho(\partial_{t}\partial_{h}u+u\cdot\nabla_{h}\partial_{h}u+w\cdot\partial_{z}\partial_{h}u+\frac{1}{Ro}\partial_{h}u^{\bot})+\rho(\partial_{h}u\cdot\nabla_{h}u+\partial_{h}w\cdot\partial_{z}u)\\
~~~~~~+\frac{\partial_{h}\rho}{\rho}\cdot(\mathbb{S}u-\frac{1}{\gamma\delta^{2}}\nabla_{h}\rho^{\gamma})+\frac{1}{\gamma\delta^{2}}\nabla_{h}\partial_{h}\rho^{\gamma}=\mathbb{S}(\partial_{h}u).
\end{cases}
\end{equation}
Taking the $L^{2}$ inner product of \eqref{eq8134}$_{2}$ with $\partial_{h}u$ and integrating by parts we then have
\begin{equation}\label{eq8136}
\begin{split}
&\frac{d}{dt}\int_{\Omega}{\frac{1}{2}\rho|\partial_{h}u|^{2}}dxdz-\int_{\Omega}{\mathbb{S}(\partial_{h}u)\cdot\partial_{h}u}dxdz+\frac{1}{\gamma\delta^{2}}\int_{\Omega}{\nabla_{h}\partial_{h}\rho^{\gamma}\cdot \partial_{h}u}dxdz\\
&=-\int_{\Omega}{\rho(\partial_{h}u\cdot\nabla_{h}u+\partial_{h}w\cdot\partial_{z}u)\cdot\partial_{h}u}dxdz-\int_{\Omega}{\frac{\partial_{h}\rho}{\rho}\cdot(\mathbb{S}u-\frac{1}{\gamma\delta^{2}}\nabla_{h}\rho^{\gamma})\cdot\partial_{h}u}dxdz.
\end{split}
\end{equation}
The term involved with pressure in \eqref{eq8136} should be carefully treated. Indeed, we can check that
\begin{equation}\label{eq8137}
\begin{split}
&\frac{1}{\gamma\delta^{2}}\int_{\Omega}{\nabla_{h}\partial_{h}\rho^{\gamma}\cdot \partial_{h}u}dxdz\\
&=\frac{1}{\delta^{2}}\int_{\Omega}{\rho^{\gamma-2}\nabla_{h}\partial_{h}\rho\cdot \rho\partial_{h}u}dxdz+\frac{1}{\delta^{2}}\int_{\Omega}{\nabla_{h}\rho^{\gamma-1}\cdot\partial_{h}\rho\cdot \partial_{h}u}dxdz\\
&=\frac{1}{\delta^{2}}\int_{\Omega}{\nabla_{h}(\rho^{\gamma-2}\partial_{h}\rho)\cdot \rho\partial_{h}u}dxdz-\frac{1}{\delta^{2}}\int_{\Omega}{\nabla_{h}\rho^{\gamma-2}\cdot\partial_{h}\rho\cdot \rho\partial_{h}u}dxdz\\
&~~~~+\frac{1}{\delta^{2}}\int_{\Omega}{\nabla_{h}\rho^{\gamma-1}\cdot\partial_{h}\rho\cdot \partial_{h}u}dxdz\\
&=-\frac{1}{\delta^{2}}\int_{\Omega}{\rho^{\gamma-2}\partial_{h}\rho\cdot\mathrm{div}_{h}\,( \rho\partial_{h}u)}dxdz+\frac{1}{(\gamma-1)\delta^{2}}\int_{\Omega}{\nabla_{h}\rho^{\gamma-1}\cdot\partial_{h}\rho\cdot \partial_{h}u}dxdz,
\end{split}
\end{equation}
where by using \eqref{eq8134}$_{1}$ we have
\begin{equation}\label{eq8135}
\begin{split}
&-\frac{1}{\delta^{2}}\int_{\Omega}{\rho^{\gamma-2}\partial_{h}\rho\cdot\mathrm{div}_{h}\,( \rho\partial_{h}u)}dxdz\\
&=\frac{1}{\delta^{2}}\int_{\Omega}{\rho^{\gamma-2}\partial_{h}\rho\cdot(\partial_{t}\partial_{h}\rho+\mathrm{div}_{h}\,(\partial_{h}\rho\cdot u)+\partial_{z}\partial_{h}(\rho w))}dxdz\\
&=\frac{d}{dt}\int_{\Omega}{\frac{1}{2\delta^{2}}\rho^{\gamma-2}|\partial_{h}\rho|^{2}}dxdz-\frac{1}{2\delta^{2}}\int_{\Omega}{\partial_{t}\rho^{\gamma-2}\cdot|\partial_{h}\rho|^{2}}dxdz\\
&~~~~+\frac{1}{\delta^{2}}\int_{\Omega}{\rho^{\gamma-2}\partial_{h}\rho\cdot(\mathrm{div}_{h}\,(\partial_{h}\rho\cdot u)+\partial_{z}\partial_{h}(\rho w))}dxdz.
\end{split}
\end{equation}
Thus, combining  \eqref{eq8137}, \eqref{eq8135} and \eqref{eq8136} we arrive at
\begin{equation*}
\begin{split}
&\frac{d}{dt}\int_{\Omega}{\frac{1}{2}\rho|\partial_{h}u|^{2}+\frac{1}{2\delta^{2}}\rho^{\gamma-2}|\partial_{h}\rho|^{2}}dxdz-\int_{\Omega}{\mathbb{S}(\partial_{h}u)\cdot\partial_{h}u}dxdz-\frac{1}{2\delta^{2}}\int_{\Omega}{\partial_{t}\rho^{\gamma-2}\cdot|\partial_{h}\rho|^{2}}dxdz\\
&+\frac{1}{\delta^{2}}\int_{\Omega}{\rho^{\gamma-2}\partial_{h}\rho\cdot(\mathrm{div}_{h}\,(\partial_{h}\rho\cdot u)+\partial_{z}\partial_{h}(\rho w))}dxdz+\frac{1}{(\gamma-1)\delta^{2}}\int_{\Omega}{\nabla_{h}\rho^{\gamma-1}\cdot\partial_{h}\rho\cdot \partial_{h}u}dxdz\\
&=-\int_{\Omega}{\rho(\partial_{h}u\cdot\nabla_{h}u+\partial_{h}w\cdot\partial_{z}u)\cdot\partial_{h}u}dxdz-\int_{\Omega}{\frac{\partial_{h}\rho}{\rho}\cdot(\mathbb{S}u-\frac{1}{\gamma\delta^{2}}\nabla_{h}\rho^{\gamma})\cdot\partial_{h}u}dxdz.
\end{split}
\end{equation*}
We now estimate the terms in the identity above. First, by \eqref{e1}$_{1}$, \eqref{eq8229}, \eqref{eq8228}, \eqref{eq82225} and \eqref{eq82230} we can check that
\begin{equation}\label{eq8231}
\begin{split}
&\frac{1}{2\delta^{2}}\int_{\Omega}{\partial_{t}\rho^{\gamma-2}\cdot|\partial_{h}\rho|^{2}}dxdz=\frac{\gamma-2}{2\delta^{2}}\int_{\Omega}{\rho^{\gamma-3}\partial_{t}\rho\cdot|\partial_{h}\rho|^{2}}dxdz\\
&=-\frac{\gamma-2}{2\delta^{2}}\int_{\Omega}{\rho^{\gamma-3}\cdot(\mathrm{div}_{h}\,(\rho u)+\partial_{z}(\rho w))\cdot|\partial_{h}\rho|^{2}}dxdz\\
&\leq \frac{C}{\delta^{2}}\big(\|\partial_{h}\rho\|_{L^{4}(\Omega)}\|u\|_{L^{4}(\Omega)}+\|\mathrm{div}_{h}\,u\|_{L^{2}(\Omega)}+\|\partial_{z}\rho\|_{L^{\infty}(\Omega)}\|w\|_{L^{2}(\Omega)}\\
&~~~~~+\|\partial_{z}w\|_{L^{2}(\Omega)}\big)\|\partial_{h}\rho\|^{2}_{L^{4}(\Omega)}\\
&\leq C\big(\delta\mathcal{A}^{2}+\mathcal{A}+\mathcal{A}^{3}\big)\mathcal{A}^{2}\leq C\big(\mathcal{A}^{2}+\mathcal{A}+\mathcal{A}^{3}\big)\mathcal{A}^{2}.
\end{split}
\end{equation}
By noting that $\rho^{\gamma-2}\partial_{h}\rho=\frac{1}{\gamma-1}\partial_{h}\rho^{\gamma-1}=\frac{\delta}{\gamma-1}\partial_{h}r$ we have
\begin{equation*}
\begin{split}
&\frac{1}{\delta^{2}}\int_{\Omega}{\rho^{\gamma-2}\partial_{h}\rho\cdot\mathrm{div}_{h}\,(\partial_{h}\rho\cdot u)}dxdz+\frac{1}{(\gamma-1)\delta^{2}}\int_{\Omega}{\nabla_{h}\rho^{\gamma-1}\cdot\partial_{h}\rho\cdot \partial_{h}u}dxdz\\
&\leq \frac{C}{\delta}\|\partial_{h}r\|_{L^{4}(\Omega)}\big(\|u\|_{L^{4}(\Omega)}\|\partial_{h}^{2}\rho\|_{L^{2}(\Omega)}+\|\partial_{h}\rho\|_{L^{4}(\Omega)}\|\partial_{h}u\|_{L^{2}(\Omega)}\big)\leq C\mathcal{A}^{3}.
\end{split}
\end{equation*}
Then, thanks to the hydrostatic balance equation \eqref{e1}$_{3}$ and the boundary condition of $w$ in \eqref{bd}, we can check that
\begin{equation*}
\begin{split}
&\frac{1}{\delta^{2}}\int_{\Omega}{\rho^{\gamma-2}\partial_{h}\rho\cdot\partial_{z}\partial_{h}(\rho w)}dxdz=\frac{1}{(\gamma-1)\delta^{2}}\int_{\Omega}{\partial_{h}\rho^{\gamma-1}\cdot\partial_{z}\partial_{h}(\rho w)}dxdz\\
&=-\frac{1}{(\gamma-1)\delta^{2}}\int_{\Omega}{\partial_{h}\partial_{z}\rho^{\gamma-1}\cdot\partial_{h}(\rho w)}dxdz=0
\end{split}
\end{equation*}
by integration by parts. Finally, direct estimations combined with \eqref{eq8228}, \eqref{eq82225} and \eqref{eq82230} yield
\begin{equation*}
\begin{split}
\int_{\Omega}{\rho(\partial_{h}u\cdot\nabla_{h}u+\partial_{h}w\cdot\partial_{z}u)\cdot\partial_{h}u}dxdz\leq C\|u\|_{H^{2}(\Omega)}^{2}\big(\|u\|_{H^{1}(\Omega)}+\|w\|_{H^{1}(\Omega)}\big)\leq C\mathcal{A}^{2}\big(\mathcal{A}+\mathcal{A}^{3}\big),
\end{split}
\end{equation*}
and
\begin{equation*}
\begin{split}
\int_{\Omega}{\frac{\partial_{h}\rho}{\rho}\cdot(\mathbb{S}u-\frac{1}{\gamma\delta^{2}}\nabla_{h}\rho^{\gamma})\cdot\partial_{h}u}dxdz&\leq C\|\partial_{h}\rho\|_{L^{4}(\Omega)}\|\partial_{h}u\|_{L^{4}(\Omega)}\big(\|\mathbb{S}u\|_{L^{2}(\Omega)}+\frac{1}{\delta^{2}}\|\partial_{h}\rho^{\gamma}\|_{L^{2}(\Omega)}\big)\\
&\leq C\delta\mathcal{A}^{2}\big(\mathcal{A}+\frac{1}{\delta}\mathcal{A}\big)\leq C\mathcal{A}^{3}.
\end{split}
\end{equation*}
Collecting the estimates above we obtain with
\begin{equation}\label{eqs1}
\begin{split}
&\frac{d}{dt}\int_{\Omega}{\frac{1}{2}\rho|\partial_{h}u|^{2}+\frac{1}{2\delta^{2}}\rho^{\gamma-2}|\partial_{h}\rho|^{2}}dxdz-\int_{\Omega}{\mathbb{S}(\partial_{h}u)\cdot\partial_{h}u}dxdz\leq C\big(\mathcal{A}^{3}+\mathcal{A}^{4}+\mathcal{A}^{5}\big).
\end{split}
\end{equation}

\textbf{Estimates~for~$(\partial_{h}^{2}r,\partial_{h}^{2}u)$.} Applying $\partial_{h}^{2}$ to \eqref{e1}$_{1,2}$ we have
\begin{equation}\label{eq8151}
\begin{cases}
\partial_{t}\partial_{h}^{2}\rho+\mathrm{div}_{h}\,(\rho\partial_{h}^{2}u)+2\mathrm{div}_{h}\,(\partial_{h}\rho\cdot\partial_{h}u)+\mathrm{div}_{h}\,(\partial_{h}^{2}\rho\cdot u)+\partial_{z}\partial_{h}^{2}(\rho w)=0,\\
\rho(\partial_{t}\partial_{h}^{2}u+u\cdot\nabla_{h}\partial_{h}^{2}u+w\cdot\partial_{z}\partial_{h}^{2}u+\frac{1}{Ro}\partial_{h}^{2}u^{\bot})+2\rho(\partial_{h}u\cdot\nabla_{h}\partial_{h}u+\partial_{h}w\cdot\partial_{z}\partial_{h}u)\\
~~~~+\rho(\partial_{h}^{2}u\cdot\nabla_{h}u+\partial_{h}^{2}w\cdot\partial_{z}u)+\frac{\partial^{2}_{h}\rho}{\rho}\cdot(\mathbb{S}u-\frac{1}{\gamma\delta^{2}}\nabla_{h}\rho^{\gamma})\\
~~~~+2\partial_{h}\rho\cdot\partial_{h}\big(\frac{1}{\rho}\cdot(\mathbb{S}u-\frac{1}{\gamma\delta^{2}}\nabla_{h}\rho^{\gamma})\big)+\frac{1}{\gamma\delta^{2}}\nabla_{h}\partial_{h}^{2}\rho^{\gamma}=\mathbb{S}(\partial_{h}^{2}u).
\end{cases}
\end{equation}
Taking the $L^{2}$ inner product of \eqref{eq8151}$_{2}$ with $\partial_{h}^{2}u$, integrating by parts
and noting that, in similar fashion to \eqref{eq8137} and \eqref{eq8135},
\begin{align}
&\frac{1}{\gamma\delta^{2}}\int_{\Omega}{\nabla_{h}\partial_{h}^{2}\rho^{\gamma}\cdot\partial_{h}^{2}u}dxdz \nonumber\\
&=\frac{1}{\delta^{2}}\int_{\Omega}{\nabla_{h}(\rho^{\gamma-1}\cdot\partial_{h}^{2}\rho)\cdot\partial_{h}^{2}u}dxdz+\frac{1}{\delta^{2}}\int_{\Omega}{\nabla_{h}(\partial_{h}\rho^{\gamma-1}\cdot\partial_{h}\rho)\cdot\partial_{h}^{2}u}dxdz \nonumber\\
&=\frac{1}{\delta^{2}}\int_{\Omega}{\nabla_{h}(\rho^{\gamma-2}\cdot\partial_{h}^{2}\rho)\cdot\rho\partial_{h}^{2}u}dxdz+\frac{1}{\delta^{2}}\int_{\Omega}{\rho^{\gamma-2}\cdot\partial_{h}^{2}\rho\cdot\nabla_{h}\rho\cdot\partial_{h}^{2}u}dxdz \nonumber\\
&~~~~+\frac{1}{\delta^{2}}\int_{\Omega}{\nabla_{h}(\partial_{h}\rho^{\gamma-1}\cdot\partial_{h}\rho)\cdot\partial_{h}^{2}u}dxdz \nonumber\\
&=\frac{1}{\delta^{2}}\int_{\Omega}{\rho^{\gamma-2}\cdot\partial_{h}^{2}\rho\cdot\big(\partial_{t}\partial_{h}^{2}\rho+2\mathrm{div}_{h}\,(\partial_{h}\rho\cdot\partial_{h}u)+\mathrm{div}_{h}\,(\partial_{h}^{2}\rho\cdot u)+\partial_{z}\partial_{h}^{2}(\rho w)\big)}dxdz \nonumber\\
&~~~~+\frac{1}{\delta^{2}}\int_{\Omega}{\rho^{\gamma-2}\cdot\partial_{h}^{2}\rho\cdot\nabla_{h}\rho\cdot\partial_{h}^{2}u}dxdz+\frac{1}{\delta^{2}}\int_{\Omega}{\nabla_{h}(\partial_{h}\rho^{\gamma-1}\cdot\partial_{h}\rho)\cdot\partial_{h}^{2}u}dxdz \nonumber\\
&=\frac{d}{dt}\int_{\Omega}{\frac{1}{2\delta^{2}}\rho^{\gamma-2}|\partial_{h}^{2}\rho|^{2}}dxdz-\frac{1}{2\delta^{2}}\int_{\Omega}{\partial_{t}\rho^{\gamma-2}\cdot|\partial^{2}_{h}\rho|^{2}}dxdz \nonumber\\
&~~~~+\frac{1}{\delta^{2}}\int_{\Omega}{\rho^{\gamma-2}\cdot\partial_{h}^{2}\rho\cdot\big(2\mathrm{div}_{h}\,(\partial_{h}\rho\cdot\partial_{h}u)+\mathrm{div}_{h}\,(\partial_{h}^{2}\rho\cdot u)+\partial_{z}\partial_{h}^{2}(\rho w)\big)}dxdz \nonumber\\
&~~~~+\frac{1}{\delta^{2}}\int_{\Omega}{\rho^{\gamma-2}\cdot\partial_{h}^{2}\rho\cdot\nabla_{h}\rho\cdot\partial_{h}^{2}u}dxdz+\frac{1}{\delta^{2}}\int_{\Omega}{\nabla_{h}(\partial_{h}\rho^{\gamma-1}\cdot\partial_{h}\rho)\cdot\partial_{h}^{2}u}dxdz,\nonumber
\end{align}
we get
\begin{align*}
&\frac{d}{dt}\int_{\Omega}{\frac{1}{2}\rho|\partial_{h}^{2}u|^{2}+\frac{1}{2\delta^{2}}\rho^{\gamma-2}|\partial_{h}^{2}\rho|^{2}}dxdz-\int_{\Omega}{\mathbb{S}(\partial_{h}^{2}u)\cdot\partial_{h}^{2}u}dxdz-\frac{1}{2\delta^{2}}\int_{\Omega}{\partial_{t}\rho^{\gamma-2}\cdot|\partial_{h}^{2}\rho|^{2}}dxdz\\
&~~~+\frac{1}{\delta^{2}}\int_{\Omega}{\rho^{\gamma-2}\cdot\partial_{h}^{2}\rho\cdot\big(2\mathrm{div}_{h}\,(\partial_{h}\rho\cdot\partial_{h}u)+\mathrm{div}_{h}\,(\partial_{h}^{2}\rho\cdot u)+\partial_{z}\partial_{h}^{2}(\rho w)\big)}dxdz\\
&~~~+\frac{1}{\delta^{2}}\int_{\Omega}{\rho^{\gamma-2}\cdot\partial_{h}^{2}\rho\cdot\nabla_{h}\rho\cdot\partial_{h}^{2}u}dxdz+\frac{1}{\delta^{2}}\int_{\Omega}{\nabla_{h}(\partial_{h}\rho^{\gamma-1}\cdot\partial_{h}\rho)\cdot\partial_{h}^{2}u}dxdz\\
&=-\int_{\Omega}{2\rho(\partial_{h}u\cdot\nabla_{h}\partial_{h}u+\partial_{h}w\cdot\partial_{z}\partial_{h}u)\cdot\partial_{h}^{2}u}dxdz-\int_{\Omega}{\rho(\partial_{h}^{2}u\cdot\nabla_{h}u+\partial_{h}^{2}w\cdot\partial_{z}u)\cdot\partial_{h}^{2}u}dxdz\\
&~~~-\int_{\Omega}{\frac{\partial^{2}_{h}\rho}{\rho}\cdot(\mathbb{S}u-\frac{1}{\gamma\delta^{2}}\nabla_{h}\rho^{\gamma})\cdot\partial_{h}^{2}u}dxdz-\int_{\Omega}{2\partial_{h}\rho\cdot\partial_{h}\big(\frac{1}{\rho}\cdot(\mathbb{S}u-\frac{1}{\gamma\delta^{2}}\nabla_{h}\rho^{\gamma})\big)\cdot\partial_{h}^{2}u}dxdz.
\end{align*}
Now we estimate the terms in the identity above. First, in similar fashion to \eqref{eq8231} we have
\begin{equation}\label{eq9231}
\begin{split}
&\frac{1}{2\delta^{2}}\int_{\Omega}{\partial_{t}\rho^{\gamma-2}\cdot|\partial_{h}^{2}\rho|^{2}}dxdz=-\frac{\gamma-2}{2\delta^{2}}\int_{\Omega}{\rho^{\gamma-3}\cdot(\mathrm{div}_{h}\,(\rho u)+\partial_{z}(\rho w))\cdot|\partial_{h}^{2}\rho|^{2}}dxdz\\
&\leq \frac{C}{\delta^{2}}\big(\|\mathrm{div}_{h}\,u\|_{L^{2}(\Omega)}+\|u\|_{L^{4}(\Omega)}\|\nabla_{h}\rho\|_{L^{4}(\Omega)}+\|\partial_{z}w\|_{L^{2}(\Omega)}+\|w\|_{L^{2}(\Omega)}\big)\|\partial_{h}^{2}\rho\|_{L^{4}(\Omega)}^{2}\\
&\leq C\big(\mathcal{A}+\mathcal{A}^{2}+\mathcal{A}^{3}\big)\times\big(\mathcal{A}+\mathcal{A}^{2}\big)^{2}\leq C\big(\mathcal{A}^{3}+\mathcal{A}^{7}\big).
\end{split}
\end{equation}
From \eqref{eq8228} and \eqref{eq82225} we obtain
\begin{equation*}
\begin{split}
&\frac{1}{\delta^{2}}\int_{\Omega}{\rho^{\gamma-2}\cdot\partial_{h}^{2}\rho\cdot\big(2\mathrm{div}_{h}\,(\partial_{h}\rho\cdot\partial_{h}u)+\mathrm{div}_{h}\,(\partial_{h}^{2}\rho\cdot u)\big)}dxdz\\
&\leq \frac{C}{\delta^{2}}\|\partial_{h}^{2}\rho\|_{L^{2}(\Omega)}\big(\|\partial_{h}\rho\|_{L^{4}(\Omega)}\|\partial_{h}^{2}u\|_{L^{4}(\Omega)}+\|\partial_{h}^{2}\rho\|_{L^{2}(\Omega)}\|\partial_{h}u\|_{L^{\infty}(\Omega)}+\|\partial_{h}^{3}\rho\|_{L^{2}(\Omega)}\|u\|_{L^{\infty}(\Omega)}\big)\\
&\leq C\mathcal{A}\big(\mathcal{A}\mathcal{B}+\mathcal{A}^{2}\big).
\end{split}
\end{equation*}
Then, by \eqref{eq8229} we get
\begin{equation*}
\begin{split}
&\frac{1}{\delta^{2}}\int_{\Omega}{\rho^{\gamma-2}\cdot\partial_{h}^{2}\rho\cdot\partial_{z}\partial_{h}^{2}(\rho w)}dxdz\\
&=\frac{C}{\delta}\int_{\Omega}{\partial_{h}^{2}r\cdot\partial_{z}\partial_{h}^{2}(\rho w)}dxdz+C\int_{\Omega}{\rho^{1-\gamma}|\partial_{h}r|^{2}\cdot\partial_{z}\partial_{h}^{2}(\rho w)}dxdz
\end{split}
\end{equation*}
for some different constants $C$ which only depend on $\gamma$. Thanks to \eqref{eq2}$_{3}$, the first term in the right side above is identically zero after integrating by parts. Thus
\begin{equation}\label{eq9232}
\begin{split}
&\frac{1}{\delta^{2}}\int_{\Omega}{\rho^{\gamma-2}\cdot\partial_{h}^{2}\rho\cdot\partial_{z}\partial_{h}^{2}(\rho w)}dxdz=-C\int_{\Omega}{\partial_{z}\rho^{1-\gamma}\cdot|\partial_{h}r|^{2}\cdot\partial_{h}^{2}(\rho w)}dxdz\\
&\leq C\|\partial_{z}\rho^{1-\gamma}\|_{L^{\infty}(\Omega)}\|\partial_{h}r\|_{L^{4}(\Omega)}^{2}\|\rho\|_{H^{2}(\Omega)}\|w\|_{H^{2}(\Omega)}\\
&\leq C\mathcal{A}^{2}\big(1+\delta\mathcal{A}\big)\times\big(\mathcal{A}+\mathcal{A}^{4}\big)\leq C\big(\mathcal{A}^{3}+\mathcal{A}^{7}\big).
\end{split}
\end{equation}
Finally, using \eqref{eq82225} and \eqref{eq8232} we get
\begin{equation*}
\begin{split}
&\frac{1}{\delta^{2}}\int_{\Omega}{\rho^{\gamma-2}\cdot\partial_{h}^{2}\rho\cdot\nabla_{h}\rho\cdot\partial_{h}^{2}u}dxdz+\frac{1}{\delta^{2}}\int_{\Omega}{\nabla_{h}(\partial_{h}\rho^{\gamma-1}\cdot\partial_{h}\rho)\cdot\partial_{h}^{2}u}dxdz\\
&\leq \frac{C}{\delta^{2}}\big(\|\partial_{h}^{2}\rho\|_{L^{2}(\Omega)}\|\partial_{h}\rho\|_{L^{4}(\Omega)}+\|\delta\partial_{h}^{2}r\|_{L^{2}(\Omega)}\|\partial_{h}\rho\|_{L^{4}(\Omega)}+\|\delta\partial_{h}r\|_{L^{4}(\Omega)}\|\partial_{h}^{2}\rho\|_{L^{2}(\Omega)}\big)\|\partial_{h}^{2}u\|_{H^{1}(\Omega)}\\
&\leq \mathcal{A}^{2}\mathcal{B},
\end{split}
\end{equation*}
\begin{equation*}
\begin{split}
&-\int_{\Omega}{2\rho(\partial_{h}u\cdot\nabla_{h}\partial_{h}u+\partial_{h}w\cdot\partial_{z}\partial_{h}u)\cdot\partial_{h}^{2}u}dxdz-\int_{\Omega}{\rho(\partial_{h}^{2}u\cdot\nabla_{h}u+\partial_{h}^{2}w\cdot\partial_{z}u)\cdot\partial_{h}^{2}u}dxdz\\
&\leq C\big(\|\partial_{h}u\|_{H^{1}(\Omega)}\|\partial_{h}^{2}u\|_{L^{2}(\Omega)}+\|\partial_{z}w\|_{H^{1}(\Omega)}\|\partial_{z}\partial_{h}u\|_{L^{2}(\Omega)}+\|\partial_{h}^{2}w\|_{L^{2}(\Omega)}\|\partial_{z}u\|_{H^{1}(\Omega)}\big)\|\partial_{h}^{2}u\|_{H^{1}(\Omega)}\\
&\leq C\big(\mathcal{A}^{2}+\big(\mathcal{A}+\mathcal{A}^{4}\big)\mathcal{A}\big)\mathcal{B}\leq C\big(\mathcal{A}^{2}+\mathcal{A}^{5}\big)\mathcal{B},
\end{split}
\end{equation*}
and
\begin{equation*}
\begin{split}
&-\int_{\Omega}{\frac{\partial^{2}_{h}\rho}{\rho}\cdot(\mathbb{S}u-\frac{1}{\gamma\delta^{2}}\nabla_{h}\rho^{\gamma})\cdot\partial_{h}^{2}u}dxdz-\int_{\Omega}{2\partial_{h}\rho\cdot\partial_{h}\big(\frac{1}{\rho}\cdot(\mathbb{S}u-\frac{1}{\gamma\delta^{2}}\nabla_{h}\rho^{\gamma})\big)\cdot\partial_{h}^{2}u}dxdz\\
&=\int_{\Omega}{\frac{\partial^{2}_{h}\rho}{\rho}\cdot(\mathbb{S}u-\frac{1}{\gamma\delta^{2}}\nabla_{h}\rho^{\gamma})\cdot\partial_{h}^{2}u}dxdz+\int_{\Omega}{\frac{2\partial_{h}\rho}{\rho}\cdot(\mathbb{S}u-\frac{1}{\gamma\delta^{2}}\nabla_{h}\rho^{\gamma})\cdot\partial_{h}^{3}u}dxdz\\
&\leq C\big(\|\partial_{h}^{2}\rho\|_{L^{2}(\Omega)}\|\partial_{h}^{2}u\|_{L^{4}(\Omega)}+\|\partial_{h}\rho\|_{L^{4}(\Omega)}\|\partial_{h}^{3}u\|_{L^{2}(\Omega)}\big)\times\big(\|\mathbb{S}u\|_{L^{4}(\Omega)}+\frac{1}{\delta^{2}}\|\nabla_{h}\rho^{\gamma}\|_{L^{4}(\Omega)}\big)\\
&\leq C\delta\mathcal{A}\mathcal{B}\big(\mathcal{B}+\frac{1}{\delta}\mathcal{A}\big)\leq C\big(\delta\mathcal{A}\mathcal{B}^{2}+\mathcal{A}^{2}\mathcal{B}\big).
\end{split}
\end{equation*}
Collecting the discussions above we conclude with
\begin{equation}\label{eqs3}
\begin{split}
&\frac{d}{dt}\int_{\Omega}{\frac{1}{2}\rho|\partial_{h}^{2}u|^{2}+\frac{1}{2\delta^{2}}\rho^{\gamma-2}|\partial_{h}^{2}\rho|^{2}}dxdz-\int_{\Omega}{\mathbb{S}(\partial_{h}^{2}u)\cdot\partial_{h}^{2}u}dxdz\\
&\leq C\big(\mathcal{A}^{3}+\mathcal{A}^{7}+\mathcal{A}^{2}\mathcal{B}+\mathcal{A}^{5}\mathcal{B}+\delta\mathcal{A}\mathcal{B}^{2}\big).
\end{split}
\end{equation}

\textbf{Estimates~for~$(\partial_{h}^{3}r,\partial_{h}^{3}u)$.} By applying $\partial_{h}^{3}$ to \eqref{e1}$_{1,2}$ we get
\begin{equation}\label{eq8171}
\begin{cases}
\partial_{t}\partial_{h}^{3}\rho+\mathrm{div}_{h}\,(\rho\partial_{h}^{3}u)+\mathrm{div}_{h}\,(3\partial_{h}\rho\cdot\partial_{h}^{2}u+3\partial_{h}^{2}\rho\cdot\partial_{h}u+\partial_{h}^{3}\rho\cdot u)+\partial_{z}\partial_{h}^{3}(\rho w)=0,\\
\rho(\partial_{t}\partial_{h}^{3}u+u\cdot\nabla_{h}\partial_{h}^{3}u+w\cdot\partial_{z}\partial_{h}^{3}u+\frac{1}{Ro}\partial_{h}^{3}u^{\bot})+3\rho(\partial_{h}u\cdot\nabla_{h}\partial_{h}^{2}u+\partial_{h}w\cdot\partial_{z}\partial_{h}^{2}u)\\
~~~~~+3\rho(\partial_{h}^{2}u\cdot\nabla_{h}\partial_{h}u+\partial_{h}^{2}w\cdot\partial_{z}\partial_{h}u)+\rho(\partial_{h}^{3}u\cdot\nabla_{h}u+\partial_{h}^{3}w\cdot\partial_{z}u)\\
~~~~~+3\partial_{h}\big[\partial_{h}\rho\cdot\partial_{h}\big(\frac{1}{\rho}(\mathbb{S}u-\frac{1}{\gamma\delta^{2}}\nabla_{h}\rho^{\gamma})\big)\big]+\partial_{h}^{3}\rho\cdot\frac{1}{\rho}(\mathbb{S}u-\frac{1}{\gamma\delta^{2}}\nabla_{h}\rho^{\gamma})\\
~~~~~+\frac{1}{\gamma\delta^{2}}\nabla_{h}\partial_{h}^{3}\rho^{\gamma}=\mathbb{S}(\partial_{h}^{3}u).
\end{cases}
\end{equation}
Taking the $L^{2}$ inner product of \eqref{eq8171}$_{2}$ with $\partial_{h}^{3}u$, integrating by parts and noting that
\begin{align*}
&\frac{1}{\gamma\delta^{2}}\int_{\Omega}{\nabla_{h}\partial_{h}^{3}\rho^{\gamma}\cdot\partial_{h}^{3}u}dxdz=\frac{1}{\delta^{2}}\int_{\Omega}{\nabla_{h}\partial_{h}^{2}(\rho^{\gamma-1}\partial_{h}\rho)\cdot\partial_{h}^{3}u}dxdz\\
&=\frac{1}{\delta^{2}}\int_{\Omega}{\rho\nabla_{h}(\rho^{\gamma-2}\cdot\partial_{h}^{3}\rho)\cdot\partial_{h}^{3}u}dxdz+\frac{1}{\delta^{2}}\int_{\Omega}{\nabla_{h}\rho\cdot\rho^{\gamma-2}\cdot\partial_{h}^{3}\rho\cdot\partial_{h}^{3}u}dxdz\\
&~~~+\frac{1}{\delta^{2}}\int_{\Omega}{\big(2\nabla_{h}(\partial_{h}\rho^{\gamma-1}\cdot\partial_{h}^{2}\rho)+\nabla_{h}(\partial_{h}^{2}\rho^{\gamma-1}\cdot\partial_{h}\rho\big)\cdot\partial_{h}^{3}u}dxdz\\
&=\frac{d}{dt}\int_{\Omega}{\frac{1}{2\delta^{2}}\rho^{\gamma-2}|\partial_{h}^{3}\rho|^{2}}dxdz-\frac{1}{2\delta^{2}}\int_{\Omega}{\partial_{t}\rho^{\gamma-2}\cdot|\partial_{h}^{3}\rho|^{2}}dxdz+\frac{1}{\delta^{2}}\int_{\Omega}{\rho^{\gamma-2}\nabla_{h}\rho\cdot\partial_{h}^{3}\rho\cdot\partial_{h}^{3}u}dxdz\\
&~~~+\frac{1}{\delta^{2}}\int_{\Omega}{\rho^{\gamma-2}\partial_{h}^{3}\rho\cdot\big(\mathrm{div}_{h}\,(3\partial_{h}\rho\cdot\partial_{h}^{2}u+3\partial_{h}^{2}\rho\cdot\partial_{h}u+\partial_{h}^{3}\rho\cdot u)+\partial_{z}\partial_{h}^{3}(\rho w)\big)}dxdz\\
&~~~+\frac{1}{\delta^{2}}\int_{\Omega}{\big(2\nabla_{h}(\partial_{h}\rho^{\gamma-1}\cdot\partial_{h}^{2}\rho)+\nabla_{h}(\partial_{h}^{2}\rho^{\gamma-1}\cdot\partial_{h}\rho\big)\cdot\partial_{h}^{3}u}dxdz,
\end{align*}
we arrive at
\begin{align*}
&\frac{d}{dt}\int_{\Omega}{\frac{1}{2}\rho|\partial_{h}^{3}u|^{2}+\frac{1}{2\delta^{2}}\rho^{\gamma-2}|\partial_{h}^{3}\rho|^{2}}dxdz-\int_{\Omega}{\mathbb{S}(\partial_{h}^{3}u)\cdot\partial_{h}^{3}u}dxdz-\frac{1}{2\delta^{2}}\int_{\Omega}{\partial_{t}\rho^{\gamma-2}\cdot|\partial_{h}^{3}\rho|^{2}}dxdz\\
&~~~+\frac{1}{\delta^{2}}\int_{\Omega}{\rho^{\gamma-2}\partial_{h}^{3}\rho\cdot\big(\mathrm{div}_{h}\,(3\partial_{h}\rho\cdot\partial_{h}^{2}u+3\partial_{h}^{2}\rho\cdot\partial_{h}u+\partial_{h}^{3}\rho\cdot u)+\partial_{z}\partial_{h}^{3}(\rho w)\big)}dxdz\\
&~~~+\frac{1}{\delta^{2}}\int_{\Omega}{\rho^{\gamma-2}\nabla_{h}\rho\cdot\partial_{h}^{3}\rho\cdot\partial_{h}^{3}u}dxdz+\frac{1}{\delta^{2}}\int_{\Omega}{\big(2\nabla_{h}(\partial_{h}\rho^{\gamma-1}\cdot\partial_{h}^{2}\rho)+\nabla_{h}(\partial_{h}^{2}\rho^{\gamma-1}\cdot\partial_{h}\rho\big)\cdot\partial_{h}^{3}u}dxdz\\
&=-\int_{\Omega}{3\rho(\partial_{h}u\cdot\nabla_{h}\partial_{h}^{2}u+\partial_{h}w\cdot\partial_{z}\partial_{h}^{2}u)\cdot\partial_{h}^{3}u}dxdz-\int_{\Omega}{\rho(\partial_{h}^{3}u\cdot\nabla_{h}u+\partial_{h}^{3}w\cdot\partial_{z}u)\cdot\partial_{h}^{3}u}dxdz\\
&~~~-\int_{\Omega}{3\rho(\partial_{h}^{2}u\cdot\nabla_{h}\partial_{h}u+\partial_{h}^{2}w\cdot\partial_{z}\partial_{h}u)\cdot\partial_{h}^{3}u}dxdz-\int_{\Omega}{\partial_{h}^{3}\rho\cdot\frac{1}{\rho}(\mathbb{S}u-\frac{1}{\gamma\delta^{2}}\nabla_{h}\rho^{\gamma})\cdot\partial_{h}^{3}u}dxdz\\
&~~~-\int_{\Omega}{3\partial_{h}\big[\partial_{h}\rho\cdot\partial_{h}\big(\frac{1}{\rho}(\mathbb{S}u-\frac{1}{\gamma\delta^{2}}\nabla_{h}\rho^{\gamma})\big)\big]\cdot\partial_{h}^{3}u}dxdz.
\end{align*}
Now, by using integration by parts, \eqref{eq82225} and \eqref{eq8232} we can verify that
\begin{equation*}
\begin{split}
&\frac{1}{2\delta^{2}}\int_{\Omega}{\partial_{t}\rho^{\gamma-2}\cdot|\partial_{h}^{3}\rho|^{2}}dxdz=-\frac{\gamma-2}{2\delta^{2}}\int_{\Omega}{\rho^{\gamma-3}\cdot(\mathrm{div}_{h}\,(\rho u)+\partial_{z}(\rho w))\cdot|\partial_{h}^{3}\rho|^{2}}dxdz\\
&\leq\frac{C}{\delta^{2}}\Big(\big(\|u\|_{H^{2}(\Omega)}\|\nabla_{h}\rho\|_{L^{\infty}(\Omega)}+\|\mathrm{div}_{h}\,u\|_{L^{\infty}(\Omega)}+\|\partial_{z}\rho^{\gamma-3}\|_{L^{\infty}(\Omega)}\|w\|_{H^{2}(\Omega)}\big)\|\partial_{h}^{3}\rho\|^{2}_{L^{2}(\Omega)}\\
&~~~~+\|w\|_{H^{2}(\Omega)}\|\partial_{h}^{3}\rho\|_{L^{2}(\Omega)}\|\partial_{z}\partial_{h}^{3}\rho\|_{L^{2}(\Omega)}\Big)\\
&\leq C\Big(\big(\delta\mathcal{A}^{2}+\mathcal{B}+\mathcal{A}+\mathcal{A}^{4}\big)\mathcal{A}^{2}+\big(\mathcal{A}+\mathcal{A}^{4}\big)\times\mathcal{A}\times\big(\mathcal{A}+\mathcal{A}^{2}+\mathcal{A}^{3}\big)\Big)\\
&\leq C\big(\mathcal{A}^{3}+\mathcal{A}^{8}+\mathcal{A}^{2}\mathcal{B}\big),
\end{split}
\end{equation*}
\begin{equation*}
\begin{split}
&\frac{1}{\delta^{2}}\int_{\Omega}{\rho^{\gamma-2}\partial_{h}^{3}\rho\cdot\mathrm{div}_{h}\,(3\partial_{h}\rho\cdot\partial_{h}^{2}u+3\partial_{h}^{2}\rho\cdot\partial_{h}u)}dxdz\\
&\leq \frac{C}{\delta^{2}}\|\partial_{h}^{3}\rho\|_{L^{2}(\Omega)}\big(\|\partial_{h}\rho\|_{L^{\infty}(\Omega)}\|\partial_{h}^{3}u\|_{L^{2}(\Omega)}+\|\partial_{h}^{2}\rho\|_{L^{4}(\Omega)}\|\partial_{h}^{2}u\|_{L^{4}(\Omega)}+\|\partial_{h}^{3}\rho\|_{L^{2}(\Omega)}\|\partial_{h}u\|_{L^{\infty}(\Omega)}\big)\\
&\leq C\mathcal{A}\big(\mathcal{A}\mathcal{B}+(\mathcal{A}+\mathcal{A}^{2})\mathcal{B}+\mathcal{A}\mathcal{B}\big)\leq C\big(\mathcal{A}^{2}\mathcal{B}+\mathcal{A}^{3}\mathcal{B}\big),
\end{split}
\end{equation*}
\begin{align*}
&\frac{1}{\delta^{2}}\int_{\Omega}{\rho^{\gamma-2}\partial_{h}^{3}\rho\cdot\mathrm{div}_{h}\,(\partial_{h}^{3}\rho\cdot u)}dxdz\\
&=-\frac{1}{2\delta^{2}}\int_{\Omega}{u\cdot\nabla_{h}\rho^{\gamma-2}\cdot|\partial_{h}^{3}\rho|^{2}}dxdz+\frac{1}{2\delta^{2}}\int_{\Omega}{\rho^{\gamma-2}\cdot\mathrm{div}_{h}\,u\cdot|\partial_{h}^{3}\rho|^{2}}dxdz\\
&\leq \frac{C}{\delta^{2}}\big(\|u\|_{L^{\infty}(\Omega)}\|\nabla_{h}\rho^{\gamma-2}\|_{L^{\infty}(\Omega)}+\|\mathrm{div}_{h}\,u\|_{L^{\infty}(\Omega)}\big)\|\partial_{h}^{3}\rho\|_{L^{2}(\Omega)}^{2}\\
&\leq C\big(\delta\mathcal{A}^{2}+\mathcal{B})\mathcal{A}^{2}\leq C(\mathcal{A}^{4}+\mathcal{A}^{2}\mathcal{B}),
\end{align*}
\begin{equation*}
\begin{split}
\frac{1}{\delta^{2}}\int_{\Omega}{\rho^{\gamma-2}\nabla_{h}\rho\cdot\partial_{h}^{3}\rho\cdot\partial_{h}^{3}u}dxdz\leq \frac{C}{\delta^{2}}\|\partial_{h}\rho\|_{L^{\infty}(\Omega)}\|\partial_{h}^{3}\rho\|_{L^{2}(\Omega)}\|\partial_{h}^{3}u\|_{L^{2}(\Omega)}\leq C\mathcal{A}^{2}\mathcal{B},
\end{split}
\end{equation*}
and
\begin{align*}
&\frac{1}{\delta^{2}}\int_{\Omega}{\big(2\nabla_{h}(\partial_{h}\rho^{\gamma-1}\cdot\partial_{h}^{2}\rho)+\nabla_{h}(\partial_{h}^{2}\rho^{\gamma-1}\cdot\partial_{h}\rho\big)\cdot\partial_{h}^{3}u}dxdz\\
&\leq \frac{C}{\delta^{2}}\Big(\|\partial_{h}^{2}\rho^{\gamma-1}\|_{L^{4}(\Omega)}\|\partial_{h}^{2}\rho\|_{L^{4}(\Omega)}+\|\partial_{h}\rho^{\gamma-1}\|_{L^{\infty}(\Omega)}\|\partial_{h}^{3}\rho\|_{L^{2}(\Omega)}\\
&~~~~+\|\partial_{h}^{3}\rho^{\gamma-1}\|_{L^{2}(\Omega)}\|\partial_{h}\rho\|_{L^{\infty}(\Omega)}\Big)\|\partial_{h}^{3}u\|_{L^{2}(\Omega)}\\
&\leq C\big(\mathcal{A}^{2}+\mathcal{A}^{3}\big)\mathcal{B}.
\end{align*}
Then, by \eqref{eq8229} we get
\begin{equation*}
\begin{split}
&\frac{1}{\delta^{2}}\int_{\Omega}{\rho^{\gamma-2}\partial_{h}^{3}\rho\cdot\partial_{z}\partial_{h}^{3}(\rho w)}dxdz\\
&=\frac{C}{\delta}\int_{\Omega}{\partial_{h}^{3}r\cdot\partial_{z}\partial_{h}^{3}(\rho w)}dxdz+\int_{\Omega}{(C\cdot\rho^{1-\gamma}\partial_{h}r\cdot\partial_{h}^{2}r+\delta C\cdot\rho^{2-2\gamma}\cdot|\partial_{h}r|^{3})\cdot\partial_{z}\partial_{h}^{3}(\rho w)}dxdz
\end{split}
\end{equation*}
for some different constants $C$ which only depend on $\gamma$. Noting that the first term in the right side above disappears after integrating by parts, we have
\begin{equation*}
\begin{split}
&\frac{1}{\delta^{2}}\int_{\Omega}{\rho^{\gamma-2}\partial_{h}^{3}\rho\cdot\partial_{z}\partial_{h}^{3}(\rho w)}dxdz\\
&=\int_{\Omega}{\partial_{h}(C\cdot\partial_{z}\rho^{1-\gamma}\cdot\partial_{h}r\cdot\partial_{h}^{2}r+\delta C\cdot\partial_{z}\rho^{2-2\gamma}\cdot|\partial_{h}r|^{3})\cdot\partial_{h}^{2}(\rho w)}dxdz\\
&\leq C\Big(\delta\|\partial_{h}r\|^{2}_{L^{6}(\mathbb{T}^{2})}\|\partial_{h}^{2}r\|_{L^{6}(\mathbb{T}^{2})}+\|\partial_{h}^{2}r\|^{2}_{L^{4}(\mathbb{T}^{2})}+\|\partial_{h}r\|_{L^{\infty}(\mathbb{T}^{2})}\|\partial_{h}^{3}r\|_{L^{2}(\mathbb{T}^{2})}\\
&~~~~+\delta^{2}\|\partial_{h}r\|^{4}_{L^{8}(\mathbb{T}^{2})}\Big)\|\rho w\|_{H^{2}(\Omega)}\\
&\leq C\big(\|r\|_{H^{3}(\mathbb{T}^{2})}^{2}+\delta\|r\|^{3}_{H^{3}(\mathbb{T}^{2})}+\delta^{2}\|r\|^{4}_{H^{3}(\mathbb{T}^{2})}\big)\|\rho\|_{H^{2}(\Omega)}\|w\|_{H^{2}(\Omega)}\\
&\leq C\big(\mathcal{A}^{2}+\delta\mathcal{A}^{3}+\delta^{2}\mathcal{A}^{4}\big)\times\big(1+\delta\mathcal{A}\big)\times\big(\mathcal{A}+\mathcal{A}^{4}\big)\leq C\big(\mathcal{A}^{3}+\mathcal{A}^{9}\big).
\end{split}
\end{equation*}
Next, from \eqref{eq971} we get
\begin{equation}\label{eq981}
\begin{split}
&\|\partial_{h}^{3}u\|_{L^{3}(\Omega)}\leq C\|\partial_{h}^{3}u\|_{L^{2}(\Omega)}^{1/2}\|\partial_{h}^{3}u\|_{H^{1}(\Omega)}^{1/2}\leq C\mathcal{A}^{1/2}\mathcal{B}^{1/2},~~\mathrm{and}\\
&\|\partial_{z}u\|_{L^{\infty}(\Omega)}\leq C\|\partial_{z}u\|_{L^{2}(\Omega)}^{1/4}\|\partial_{z}u\|_{H^{2}(\Omega)}^{3/4}\leq C\mathcal{A}^{1/4}\mathcal{B}^{3/4},
\end{split}
\end{equation}
which imply that
\begin{equation*}
\begin{split}
&\int_{\Omega}{3\rho(\partial_{h}u\cdot\nabla_{h}\partial_{h}^{2}u+\partial_{h}w\cdot\partial_{z}\partial_{h}^{2}u)\cdot\partial_{h}^{3}u}dxdz\\
&\leq C\big(\|u\|_{H^{2}(\Omega)}\|\partial_{h}^{3}u\|_{L^{2}(\Omega)}+\|w\|_{H^{2}(\Omega)}\|\partial_{z}\partial_{h}^{2}u\|_{L^{2}(\Omega)}\big)\|\partial_{h}^{3}u\|_{L^{3}(\Omega)}\\
&\leq C\big(\mathcal{A}\mathcal{B}+\big(\mathcal{A}+\mathcal{A}^{4}\big)\mathcal{B}\big)\mathcal{A}^{1/2}\mathcal{B}^{1/2}\leq C\big(\mathcal{A}^{3/2}\mathcal{B}^{3/2}+\mathcal{A}^{9/2}\mathcal{B}^{3/2}\big),
\end{split}
\end{equation*}
\begin{align*}
&\int_{\Omega}{\rho(\partial_{h}^{3}u\cdot\nabla_{h}u+\partial_{h}^{3}w\cdot\partial_{z}u)\cdot\partial_{h}^{3}u}dxdz\\
&=\int_{\Omega}{\rho\partial_{h}^{3}u\cdot\nabla_{h}u\cdot\partial_{h}^{3}u}dxdz\int_{\Omega}{\rho\partial_{h}^{2}w\cdot\partial_{z}u\cdot\partial_{h}^{4}u}dxdz-\int_{\Omega}{\rho\partial_{h}^{2}w\cdot\partial_{z}\partial_{h}u\cdot\partial_{h}^{3}u}dxdz\\
&~~~~-\int_{\Omega}{\partial_{h}\rho\cdot\partial_{h}^{2}w\cdot\partial_{z}u\cdot\partial_{h}^{3}u}dxdz\\
&\leq C\|\partial_{h}^{3}u\|_{L^{2}(\Omega)}\|\nabla_{h}u\|_{L^{6}(\Omega)}\|\partial_{h}^{3}u\|_{L^{3}(\Omega)}+C\|w\|_{H^{2}(\Omega)}\|\partial_{z}u\|_{L^{\infty}(\Omega)}\|\partial_{h}^{4}u\|_{L^{2}(\Omega)}\\
&~~~~+\big(\|\partial_{h}\rho\|_{L^{\infty}(\Omega)}\|\partial_{z}u\|_{L^{6}(\Omega)}+\|\partial_{z}\partial_{h}u\|_{L^{6}(\Omega)}\big)\|w\|_{H^{2}(\Omega)}\|\partial_{h}^{3}u\|_{L^{3}(\Omega)}\\
&\leq C\big(\mathcal{A}^{3}+\big(\mathcal{A}+\mathcal{A}^{4}\big)\mathcal{A}^{1/4}\mathcal{B}^{7/4}\big)+C\big(\delta\mathcal{A}^{2}+\mathcal{B}\big)\times\big(\mathcal{A}+\mathcal{A}^{4}\big)\times\mathcal{A}\\
&\leq C\big(\mathcal{A}^{3}+\mathcal{A}^{7}+\mathcal{A}^{5/4}\mathcal{B}^{7/4}+\mathcal{A}^{17/4}\mathcal{B}^{7/4}+\mathcal{A}^{2}\mathcal{B}+\mathcal{A}^{5}\mathcal{B}\big),
\end{align*}
and
\begin{equation*}
\begin{split}
\int_{\Omega}{\partial_{h}^{3}\rho\cdot\frac{1}{\rho}(\mathbb{S}u-\frac{1}{\gamma\delta^{2}}\nabla_{h}\rho^{\gamma})\cdot\partial_{h}^{3}u}dxdz&\leq C\|\partial_{h}^{3}\rho\|_{L^{2}(\Omega)}\|\partial_{h}^{3}u\|_{L^{3}(\Omega)}\big(\|\mathbb{S}u\|_{L^{6}(\Omega)}+\frac{1}{\delta^{2}}\|\partial_{h}\rho^{\gamma}\|_{L^{6}(\Omega)}\big)\\
&\leq C\delta\mathcal{A}^{3/2}\mathcal{B}^{1/2}\big(\mathcal{B}+\frac{1}{\delta}\mathcal{A}\big)\leq C\big(\mathcal{A}^{3/2}\mathcal{B}^{3/2}+\mathcal{A}^{5/2}\mathcal{B}^{1/2}\big).
\end{split}
\end{equation*}
Finally, it can be easily verified that
\begin{equation*}
\begin{split}
&\int_{\Omega}{3\rho(\partial_{h}^{2}u\cdot\nabla_{h}\partial_{h}u+\partial_{h}^{2}w\cdot\partial_{z}\partial_{h}u)\cdot\partial_{h}^{3}u}dxdz\\
&\leq C\big(\|u\|_{H^{2}(\Omega)}+\|w\|_{H^{2}(\Omega)}\big)\|\partial_{h}u\|_{H^{2}(\Omega)}\|\partial_{h}^{3}u\|_{L^{3}(\Omega)}\leq C\big(\mathcal{A}+\mathcal{A}^{4}\big)\mathcal{B}\mathcal{A},
\end{split}
\end{equation*}
and
\begin{equation*}
\begin{split}
&\int_{\Omega}{3\partial_{h}\big[\partial_{h}\rho\cdot\partial_{h}\big(\frac{1}{\rho}(\mathbb{S}u-\frac{1}{\gamma\delta^{2}}\nabla_{h}\rho^{\gamma})\big)\big]\cdot\partial_{h}^{3}u}dxdz=-\int_{\Omega}{3\partial_{h}\rho\cdot\partial_{h}\big(\frac{1}{\rho}(\mathbb{S}u-\frac{1}{\gamma\delta^{2}}\nabla_{h}\rho^{\gamma})\big)\cdot\partial_{h}^{4}u}dxdz\\
&\leq C\|\partial_{h}\rho\|_{L^{\infty}(\Omega)}\|\partial_{h}^{4}u\|_{L^{2}(\Omega)}\Big(\|\partial_{h}\rho\|_{L^{\infty}(\Omega)}\|\mathbb{S}u\|_{L^{2}(\Omega)}+\|\partial_{h}\mathbb{S}u\|_{L^{2}(\Omega)}+\frac{1}{\delta^{2}}\|\partial_{h}^{2}\rho^{\gamma}\|_{L^{2}(\Omega)}\\
&~~~~~~~~~~+\frac{1}{\delta^{2}}\|\partial_{h}\rho\|_{L^{\infty}(\Omega)}\|\partial_{h}\rho^{\gamma}\|_{L^{2}(\Omega)}\Big)\\
&\leq C\delta\mathcal{A}\mathcal{B}\big(\delta\mathcal{A}^{2}+\mathcal{B}+\frac{1}{\delta}\big(\mathcal{A}+\mathcal{A}^{2}\big)+\mathcal{A}^{2}\big)\leq C\big(\mathcal{A}^{2}\mathcal{B}+\mathcal{A}^{3}\mathcal{B}+\delta\mathcal{A}\mathcal{B}^{2}\big).
\end{split}
\end{equation*}
Collecting the above discussions we obtain with
\begin{equation}\label{eqs4}
\begin{split}
&\frac{d}{dt}\int_{\Omega}{\frac{1}{2}\rho|\partial_{h}^{3}u|^{2}+\frac{1}{2\delta^{2}}\rho^{\gamma-2}|\partial_{h}^{3}\rho|^{2}}dxdz-\int_{\Omega}{\mathbb{S}(\partial_{h}^{3}u)\cdot\partial_{h}^{3}u}dxdz\\
&\leq C\Big(\mathcal{A}^{3}+\mathcal{A}^{4}+\mathcal{A}^{7}+\mathcal{A}^{8}+\mathcal{A}^{9}+\mathcal{A}^{2}\mathcal{B}+\mathcal{A}^{3}\mathcal{B}+\mathcal{A}^{5}\mathcal{B}\\
&~~~~+\mathcal{A}^{3/2}\mathcal{B}^{3/2}+\mathcal{A}^{9/2}\mathcal{B}^{3/2}+\mathcal{A}^{5/2}\mathcal{B}^{1/2}+\mathcal{A}^{5/4}\mathcal{B}^{7/4}+\mathcal{A}^{17/4}\mathcal{B}^{7/4}+\delta\mathcal{A}\mathcal{B}^{2}\Big).
\end{split}
\end{equation}

\textbf{Estimates~for~$(\partial_{z}u,\partial_{z}^{2}u,\partial_{h}\partial_{z}u)$.} To derive the estimates for the vertical derivatives of $u$, it's more convenient to work directly with the velocity equation in \eqref{eq2}. By applying $\partial_{z}$ to \eqref{eq2}$_{2}$ we get
\begin{equation}\label{eq8173}
\partial_{t}\partial_{z}u+\partial_{z}(u\cdot\nabla_{h} u)+\partial_{z}(w\cdot\partial_{z}u)+\frac{1}{Ro}\partial_{z}u^{\bot}-\mathbb{S}(\partial_{z}u)=\partial_{z}\Big(\big(\frac{1}{\rho}-1\big)\mathbb{S}u\Big).
\end{equation}
Taking the $L^{2}$ inner product of \eqref{eq8173} with $\partial_{z}u-\partial_{z}^{3}u$ and integrating by parts we arrive at
\begin{equation}\label{eq8241}
\begin{split}
&\frac{1}{2}\frac{d}{dt}\sum\limits_{n=1}^{2}\int_{\Omega}{|\partial_{z}^{n}u|^{2}}dxdz-\int_{\Omega}{\mathbb{S}(\partial_{z}u)\cdot\big(\partial_{z}u-\partial_{z}^{3}u\big)}dxdz\\
&=-\int_{\Omega}{\partial_{z}(u\cdot\nabla_{h}u+w\cdot\partial_{z}u)\cdot(\partial_{z}u-\partial_{z}^{3}u)}dxdz+\int_{\Omega}{\partial_{z}\Big(\big(\frac{1}{\rho}-1\big)\mathbb{S}u\Big)\cdot(\partial_{z}u-\partial_{z}^{3}u)}dxdz.
\end{split}
\end{equation}
First, by using \eqref{eq82230} and \eqref{eq981} we can verify that
\begin{align}
&\int_{\Omega}{\partial_{z}(u\cdot\nabla_{h}u+w\cdot\partial_{z}u)\cdot(\partial_{z}u-\partial_{z}^{3}u)}dxdz \nonumber\\
&\leq C\Big(\|\partial_{z}u\|_{L^{6}(\Omega)}\|\nabla_{h}u\|_{L^{3}(\Omega)}+\|u\|_{L^{\infty}(\Omega)}\|\nabla_{h}\partial_{z}u\|_{L^{2}(\Omega)}\nonumber\\
&~~~~~~~~~~+\|\partial_{z}w\|_{L^{2}(\Omega)}\|\partial_{z}u\|_{L^{\infty}(\Omega)}+\|w\|_{L^{6}(\Omega)}\|\partial^{2}_{z}u\|_{L^{3}(\Omega)}\Big)\times\|(\partial_{z}u,\partial_{z}^{3}u)\|_{L^{2}(\Omega)}\nonumber\\
&\leq C\Big(\mathcal{A}^{2}+\|w\|_{H^{1}(\Omega)}\big(\|\partial_{z}u\|^{1/4}_{L^{2}(\Omega)}\|\partial_{z}u\|_{H^{2}(\Omega)}^{3/4}+\|\partial_{z}^{2}u\|^{1/2}_{L^{2}(\Omega)}\|\partial^{2}_{z}u\|_{H^{1}(\Omega)}^{1/2}\big)\Big)\mathcal{B}\nonumber\\
&\leq C\Big(\mathcal{A}^{2}+\big(\mathcal{A}+\mathcal{A}^{3}\big)\times\big(\mathcal{A}^{1/4}\mathcal{B}^{3/4}+\mathcal{A}^{1/2}\mathcal{B}^{1/2}\big)\Big)\mathcal{B}\nonumber\\
&\leq C\big(\mathcal{A}^{2}\mathcal{B}+\mathcal{A}^{5/4}\mathcal{B}^{7/4}+\mathcal{A}^{3/2}\mathcal{B}^{3/2}+\mathcal{A}^{13/4}\mathcal{B}^{7/4}+\mathcal{A}^{7/2}\mathcal{B}^{3/2}\big),\nonumber
\end{align}
and
\begin{equation*}
\begin{split}
&\int_{\Omega}{\partial_{z}\Big(\big(\frac{1}{\rho}-1\big)\mathbb{S}u\Big)\cdot\partial_{z}u}dxdz=-\int_{\Omega}{\big(\frac{1}{\rho}-1\big)\mathbb{S}u\cdot\partial_{z}^{2}u}dxdz\leq C\|\mathbb{S}u\|_{L^{2}(\Omega)}\|\partial_{z}^{2}u\|_{L^{2}(\Omega)}\leq C\mathcal{A}^{2}.
\end{split}
\end{equation*}
To estimate the  last term in \eqref{eq8241}, we split
\begin{equation}\label{eq8183}
\begin{split}
&-\int_{\Omega}{\partial_{z}\Big(\big(\frac{1}{\rho}-1\big)\mathbb{S}u\Big)\cdot\partial_{z}^{3}u}dxdz\\
&=\int_{\Omega}{\frac{\partial_{z}\rho}{\rho^{2}}\cdot\mathbb{S}u\cdot\partial_{z}^{3}u}+\int_{\Omega}{\big(\frac{\rho-(\xi^{L})^{\frac{1}{\gamma-1}}}{\rho}\big)\cdot\partial_{z}\mathbb{S}u\cdot\partial_{z}^{3}u}+\int_{\Omega}{\big(\frac{(\xi^{L})^{\frac{1}{\gamma-1}}-1}{\rho}\big)\cdot\partial_{z}\mathbb{S}u\cdot\partial_{z}^{3}u}\\
&\triangleq J_{1}+J_{2}+J_{3},
\end{split}
\end{equation}
where we have
\begin{equation}\label{eq8182}
\begin{split}
J_{1}\leq C\|\partial_{z}\rho\|_{L^{\infty}(\Omega)}\|\mathbb{S}u\|_{L^{2}(\Omega)}\|\partial_{z}^{3}u\|_{L^{2}(\Omega)}\leq C\mathcal{A}\mathcal{B}.
\end{split}
\end{equation}
Furthermore, by using differential mean value theorem, we can check that
\begin{equation*}
\|\rho-(\xi^{L})^{\frac{1}{\gamma-1}}\|_{L^{\infty}(\Omega)}\leq C\delta \|r\|_{L^{\infty}(\Omega)}.
\end{equation*}
Thus
\begin{equation}\label{eq8181}
J_{2}\leq C\delta\|r\|_{L^{\infty}(\Omega)}\|\partial_{z}\mathbb{S}u\|_{L^{2}(\Omega)}\|\partial_{z}^{3}u\|_{L^{2}(\Omega)}\leq C\delta\mathcal{A}\mathcal{B}^{2}.
\end{equation}
Now it suffices to estimate $J_{3}$. Noting that $(\xi^{L})^{\frac{1}{\gamma-1}}-1<0$, we have, after dropping the negative term,
\begin{equation*}
\begin{split}
J_{3}\leq\int_{\Omega}{\big(\frac{(\xi^{L})^{\frac{1}{\gamma-1}}-1}{\rho}\big)\cdot\mu\triangle_{h}\partial_{z}u\cdot\partial_{z}^{3}u}+\int_{\Omega}{\big(\frac{(\xi^{L})^{\frac{1}{\gamma-1}}-1}{\rho}\big)\cdot\nu\nabla_{h}\mathrm{div}_{h}\,\partial_{z}u\cdot\partial_{z}^{3}u}\triangleq J_{31}+J_{32}.
\end{split}
\end{equation*}
Using integration by parts we get
\begin{align*}
J_{31}&=-\mu\int_{\Omega}{\frac{(\xi^{L})^{\frac{1}{\gamma-1}}-1}{\rho}\cdot\nabla_{h}\partial_{z}u\cdot\nabla_{h}\partial_{z}^{3}u}+\mu\int_{\Omega}{\frac{(\xi^{L})^{\frac{1}{\gamma-1}}-1}{\rho^{2}}\cdot\nabla_{h}\rho\cdot\nabla_{h}\partial_{z}u\cdot\partial_{z}^{3}u}\\
&\leq\mu\int_{\Omega}{\partial_{z}\big(\frac{(\xi^{L})^{\frac{1}{\gamma-1}}-1}{\rho}\big)\cdot\nabla_{h}\partial_{z}u\cdot\nabla_{h}\partial_{z}^{2}u}+\mu\int_{\Omega}{\frac{(\xi^{L})^{\frac{1}{\gamma-1}}-1}{\rho^{2}}\cdot\nabla_{h}\rho\cdot\nabla_{h}\partial_{z}u\cdot\partial_{z}^{3}u}\\
&\leq C\big\|\partial_{z}\big(\frac{(\xi^{L})^{\frac{1}{\gamma-1}}-1}{\rho}\big)\big\|_{L^{\infty}(\Omega)}\|\nabla_{h}\partial_{z}u\|_{L^{2}(\Omega)}\|\nabla_{h}\partial_{z}^{2}u\|_{L^{2}(\Omega)}\\
&~~~~+C\|\nabla_{h}\rho\|_{L^{3}(\Omega)}\|\nabla_{h}\partial_{z}u\|_{L^{6}(\Omega)}\|\partial_{z}^{3}u\|_{L^{2}(\Omega)}\\
&\leq C\big(\mathcal{A}\mathcal{B}+\delta\mathcal{A}\mathcal{B}^{2}\big).
\end{align*}
The same estimate holds for $J_{32}$. Thus we have $J_{3}\leq C\big(\mathcal{A}\mathcal{B}+\delta\mathcal{A}\mathcal{B}^{2}\big)$. Combining this estimate with \eqref{eq8183}, \eqref{eq8182} and \eqref{eq8181} we obtain
\begin{equation*}
-\int_{\Omega}{\partial_{z}\Big(\big(\frac{1}{\rho}-1\big)\mathbb{S}u\Big)\cdot\partial_{z}^{3}u}dxdz\leq C\big(\mathcal{A}\mathcal{B}+\delta\mathcal{A}\mathcal{B}^{2}\big).
\end{equation*}
Thus, collecting the estimates above we conclude with
\begin{equation}\label{eqs5}
\begin{split}
&\frac{1}{2}\frac{d}{dt}\sum\limits_{n=1}^{2}\int_{\Omega}{|\partial_{z}^{n}u|^{2}}dxdz-\int_{\Omega}{\mathbb{S}(\partial_{z}u)\cdot\big(\partial_{z}u-\partial_{z}^{3}u\big)}dxdz\\
&\leq C\big(\mathcal{A}^{2}+\mathcal{A}^{2}\mathcal{B}+\mathcal{A}\mathcal{B}+\mathcal{A}^{5/4}\mathcal{B}^{7/4}+\mathcal{A}^{3/2}\mathcal{B}^{3/2}+\mathcal{A}^{13/4}\mathcal{B}^{7/4}+\mathcal{A}^{7/2}\mathcal{B}^{3/2}+\delta\mathcal{A}\mathcal{B}^{2}\big).
\end{split}
\end{equation}

Finally, applying $\partial_{h}\partial_{z}$ to \eqref{eq2}$_{2}$, taking the $L^{2}$ inner product of the result with $\partial_{h}\partial_{z}u$ and using a similar approach as the previous discussions we can also get
\begin{equation}\label{eqs6}
\begin{split}
&\frac{1}{2}\frac{d}{dt}\int_{\Omega}{|\partial_{h}\partial_{z}u|^{2}}dxdz-\int_{\Omega}{\mathbb{S}(\partial_{h}\partial_{z}u)\cdot\partial_{h}\partial_{z}u}dxdz\\
&\leq C\big(\mathcal{A}^{2}+\mathcal{A}^{2}\mathcal{B}+\mathcal{A}\mathcal{B}+\mathcal{A}^{5/4}\mathcal{B}^{7/4}+\mathcal{A}^{3/2}\mathcal{B}^{3/2}+\mathcal{A}^{13/4}\mathcal{B}^{7/4}+\mathcal{A}^{7/2}\mathcal{B}^{3/2}+\delta\mathcal{A}\mathcal{B}^{2}\big).
\end{split}
\end{equation}

\textbf{Uniform~estimates.} Combining \eqref{eqs2}, \eqref{eqs1}, \eqref{eqs3}, \eqref{eqs4}, \eqref{eqs5} and \eqref{eqs6}, and using Young's inequalities we have
\begin{equation*}
\begin{split}
&\frac{1}{2}\frac{d}{dt}\mathcal{A}^{2}+\mathcal{B}^{2}\leq C\sum\limits_{n=2}^{9}\mathcal{A}^{n}+C\big(\mathcal{A}+\mathcal{A}^{2}+\mathcal{A}^{3}+\mathcal{A}^{5}\big)\mathcal{B}+C\mathcal{A}^{5/2}\mathcal{B}^{1/2}\\
&~~~~+C\big(\mathcal{A}^{5/4}+\mathcal{A}^{13/4}+\mathcal{A}^{17/4}\big)\mathcal{B}^{7/4}+C\big(\mathcal{A}^{3/2}+\mathcal{A}^{7/2}+\mathcal{A}^{9/2}\big)\mathcal{B}^{3/2}+C\delta\mathcal{A}\mathcal{B}^{2}\\
&\leq C\big(\sum\limits_{n=2}^{9}\mathcal{A}^{n}+\mathcal{A}^{10}+\mathcal{A}^{10/3}+\mathcal{A}^{26}+\mathcal{A}^{34}+\mathcal{A}^{14}+\mathcal{A}^{18}\big)+\frac{1}{2}\mathcal{B}^{2}+C\delta\mathcal{A}\mathcal{B}^{2}\\
&\leq C\big(\mathcal{A}^{2}+\mathcal{A}^{34}\big)+\frac{1}{2}\mathcal{B}^{2}+C\delta\mathcal{A}\mathcal{B}^{2}.
\end{split}
\end{equation*}
As a result, by \eqref{eq82210} we get
\begin{equation*}
\begin{split}
\frac{d}{dt}\mathcal{A}^{2}(t)+\mathcal{B}^{2}(t)&\leq M_{1}\big(1+\mathcal{A}^{32}\big)\mathcal{A}^{2}(t)+\delta M_{2}\mathcal{A}\mathcal{B}^{2}(t)\\
&\leq M_{1}\big(1+\mathcal{A}^{32}\big)\mathcal{A}^{2}(t)+\delta M_{2}C_{6}\mathcal{A}_{*}\mathcal{B}^{2}(t),~~\forall\,t\in (0,T_{\delta}^{*}]
\end{split}
\end{equation*}
for some constants $M_{1}$ and $M_{2}$ which are independent of $\delta$. Now we let
\begin{equation*}
\delta_{0}:=\min\Big\{\delta_{**},\frac{1}{4M_{2}C_{5}C_{6}^{2}B_{in}}\Big\},
\end{equation*}
where $\delta_{**}$ is defined by \eqref{eq9253}, $C_{5}$ and $C_{6}$ are from \eqref{eq82210}, and $B_{in}$ is from \eqref{eq974}. Then when $\delta\leq \delta_{0}$ and $t\leq T_{\delta}^{*}$ we have
\begin{equation*}
\delta M_{2}C_{6}\mathcal{A}_{*}\leq  \frac{1}{4M_{2}C_{5}C_{6}^{2}B_{in}}\times M_{2}C_{6}\times2C_{5}C_{6}B_{in}=\frac{1}{2}.
\end{equation*}
Thus
\begin{equation*}
\begin{split}
\frac{d}{dt}\mathcal{A}^{2}(t)+\frac{1}{2}\mathcal{B}^{2}(t)&\leq M_{1}\big(1+\mathcal{A}^{32}\big)\mathcal{A}^{2}(t),~~\forall\,t\in (0,T_{\delta}^{*}].
\end{split}
\end{equation*}
By applying \eqref{eq82210}, \eqref{eq8261}, Gronwall's inequality and the fact $\mathcal{A}_{*}(0)\leq B_{in}$ we get
\begin{equation*}
\begin{split}
&\frac{1}{C_{5}^{2}}\mathcal{A}_{*}^{2}(t)+\frac{1}{2}\int^{t}_{0}{\mathcal{B}^{2}(s)}ds\leq
\mathcal{A}^{2}(t)+\frac{1}{2}\int^{t}_{0}{\mathcal{B}^{2}(s)}ds\leq e^{M_{1}\int_{0}^{t}{(1+\mathcal{A}^{32})}ds}\mathcal{A}^{2}(0)\\
&\leq C_{6}^{2}e^{M_{1}\int_{0}^{t}{(1+(C_{6}\mathcal{A}_{*})^{32})}ds}\mathcal{A}_{*}^{2}(0)\leq C_{6}^{2}e^{M_{1}(1+(2C_{5}C_{6}^{2}B_{in})^{32})t}\mathcal{A}_{*}^{2}(0),~~\forall\,t\in (0,T_{\delta}^{*}].
\end{split}
\end{equation*}
Thus
\begin{equation*}
\begin{split}
&\sup\limits_{s\in[0,t]}\mathcal{A}_{*}^{2}(s)+\frac{C_{5}^{2}}{2}\int^{t}_{0}{\mathcal{B}^{2}(s)}ds\leq 2C_{5}^{2}C_{6}^{2}e^{M_{1}(1+(2C_{5}C_{6}^{2}B_{in})^{32})t}\mathcal{A}_{*}^{2}(0),~~\forall\,t\in (0,T_{\delta}^{*}].
\end{split}
\end{equation*}
Now we take $T$ such that $e^{M_{1}(1+(2C_{5}C_{6}^{2}B_{in})^{32})T}=\frac{3}{2}$. Then we obtain
\begin{equation*}
\begin{split}
&\sup\limits_{s\in[0,T]}\mathcal{A}_{*}^{2}(s)+\frac{C_{5}^{2}}{2}\int^{T}_{0}{\mathcal{B}^{2}(s)}ds\leq 3C_{5}^{2}C_{6}^{2}\mathcal{A}_{*}^{2}(0),
\end{split}
\end{equation*}
which implies \eqref{est} and $T_{\delta}^{*}\geq T$ by the definition of $T^{*}_{\delta}$ in \eqref{eq8261}. Now we have completed the proof of Theorem \ref{th1}.

\section{Low Mach number limit with well-prepared data}\label{sec3}
In this section we investigate the low Mach number singular limit of the compressible primitive equations with well-prepared initial data. Let $\{(r,u)\}_{\delta}$ be the solution of \eqref{eq2} given by Theorem \ref{th1}. From the discussions in Section \ref{sec2} we know that $\delta \|r\|_{L^{\infty}(\mathbb{T}^{2})}\leq \frac{1}{4}$. Thus, by differential mean value theorem we get
\begin{align}\label{eq943}
|\rho(t,x,z)-\rho^{L}(z)|\leq C\delta|r(t,x)|.
\end{align}
Furthermore, we need the following estimates for the solution of the limiting system \eqref{eqL}-\eqref{inL}:
\begin{equation}\label{eq9111}
\|u^{L}\|_{L^{\infty}(0,T;H^{2}(\Omega))}+\|\partial_{z}u^{L}\|_{L^{2}(0,T;H^{2}(\Omega))}\leq C,
\end{equation}
and
\begin{equation}\label{eq9112}
\|\Pi\|_{L^{\infty}(0,T;H^{2}(\mathbb{T}^{2}))}+\|\partial_{t}\Pi\|_{L^{2}(0,T;H^{1}(\mathbb{T}^{2}))}+\|\Pi|_{t=0}\|_{H^{1}(\mathbb{T}^{2})}\leq C.
\end{equation}
These estimates can be proved by adapting the discussions in Appendix in \cite{Dai-Mu} under the condition $u^{L}_{0}\in H^{3}(\Omega)$.

Now we begin to prove \eqref{eqWell}. We start from constructing the estimate for the difference $w-w^{L}$. By \eqref{eq-nota} and \eqref{eqL}$_{2}$ we can verify that
\begin{align*}
\mathrm{div}_{h}(\xi^{L}u^{L})+\partial_{z}(\xi^{L}w^{L})+\frac{\gamma-2}{\gamma-1}\theta w^{L}=0,
\end{align*}
which implies
\begin{align}\label{eq922}
w^{L}(z)=-\frac{1}{\xi^{L}}\int^{z}_{0}{\mathrm{div}_{h}(\widetilde{\xi^{L}u^{L}})}d\tau-\frac{(\gamma-2)\theta}{(\gamma-1)\xi^{L}}\int^{z}_{0}{\widetilde{w^{L}}}d\tau
\end{align}
by noting that $w^{L}|_{z=0,\kappa}=0$. Subtracting \eqref{eq922} from \eqref{eq921} we further have
\begin{align*}
&w(z)-w^{L}(z)=-\Big(\frac{1}{\xi}-\frac{1}{\xi^{L}}\Big)\int^{z}_{0}{\mathrm{div}_{h}(\widetilde{\xi u})}d\tau-\frac{1}{\xi^{L}}\int^{z}_{0}{\mathrm{div}_{h}(\widetilde{\xi u}-\widetilde{\xi^{L} u^{L}})}d\tau\\
&~~~~-\frac{(\gamma-2)\theta}{\gamma-1}\Big(\frac{1}{\xi}-\frac{1}{\xi^{L}}\Big)\int^{z}_{0}{\widetilde{w}}d\tau-\frac{(\gamma-2)\theta}{(\gamma-1)\xi^{L}}\int^{z}_{0}{\widetilde{w}-\widetilde{w^{L}}}d\tau+\frac{\delta(\gamma-2)}{(\gamma-1)\xi}\int^{z}_{0}{\widetilde{u}\cdot\nabla_{h}r}d\tau.
\end{align*}
Thus, by using \eqref{est} and \eqref{eq8253} we can obtain that
\begin{align*}
&\|w-w^{L}\|_{L^{2}(\Omega)}\leq C\|\xi-\xi^{L}\|_{L^{6}(\Omega)}\big(\|\mathrm{div}_{h}\,u\|_{L^{3}(\Omega)}+\delta\|\nabla_{h}r\|_{L^{6}(\Omega)}\|u\|_{L^{6}(\Omega)}\big)\\
&~~~+C\big(\|\mathrm{div}_{h}\,((\xi -\xi^{L})u)\|_{L^{2}(\Omega)}+\|\mathrm{div}_{h}\,(\xi^{L}( u-u^{L}))\|_{L^{2}(\Omega)}\big)+C\|\xi-\xi^{L}\|_{L^{4}(\Omega)}\|w\|_{L^{4}(\Omega)}\\
&~~~+\frac{|\gamma-2|\theta\kappa}{\sqrt{2}(\gamma-1)}\|(\xi^{L})^{-1}\|_{L^{\infty}(\Omega)}\|w-w^{L}\|_{L^{2}(\Omega)}+C\delta\|u\|_{L^{\infty}(\Omega)}\|\nabla_{h}r\|_{L^{2}(\Omega)}\\
&\leq C\delta\|r\|_{H^{1}(\Omega)}\big(\mathcal{A}+\delta\mathcal{A}^{2}\big)+C\big(\delta\|r\|_{L^{4}(\Omega)}\|\mathrm{div}_{h}\,u\|_{L^{4}}+\delta\|\nabla_{h}r\|_{L^{2}(\Omega)}\|u\|_{L^{\infty}}\big)\\
&~~~+C\|\mathrm{div}_{h}\,(u-u^{L})\|_{L^{2}(\Omega)}+C\delta\|r\|_{H^{1}(\Omega)}\|w\|_{H^{1}(\Omega)}+\frac{1}{4}\|w-w^{L}\|_{L^{2}(\Omega)}+C\delta\mathcal{A}\|r\|_{H^{1}(\Omega)}\\
&\leq C\delta\|r\|_{H^{1}(\mathbb{T}^{2})}+C\|u-u^{L}\|_{H^{1}(\Omega)}+\frac{1}{4}\|w-w^{L}\|_{L^{2}(\Omega)},
\end{align*}
which immediately leads to
\begin{align}\label{eq941}
\|w-w^{L}\|_{L^{2}(\Omega)}\leq C\big(\delta\|r\|_{H^{1}(\mathbb{T}^{2})}+\|u-u^{L}\|_{H^{1}(\Omega)}\big).
\end{align}

Now, subtracting \eqref{eqL}$_{1}$ from \eqref{eq2}$_{2}$ and denoting by
\begin{equation*}
v:=u-u^{L},
\end{equation*}
we get
\begin{equation}\label{eq914}
\begin{split}
&\partial_{t}v+u\cdot\nabla_{h}v+v\cdot\nabla_{h}u^{L}+w\cdot\partial_{z}v+(w-w^{L})\cdot\partial_{z}u^{L}+\frac{1}{Ro}v^{\bot}-\frac{1}{\rho^{L}}\mathbb{S}v\\
&=-\nabla_{h}\Big(\frac{1}{(\gamma-1)\delta}r-\Pi\Big)+\Big(\frac{1}{\rho}-\frac{1}{\rho^{L}}\Big)\mathbb{S}u.
\end{split}
\end{equation}
We introduce the notations
\begin{align*}
\mathcal{D}=\mathcal{D}(t):=\Big(\int_{\Omega}{\rho^{L}\cdot\big(|v|^{2}+|\partial_{h}v|^{2}\big)+|\partial_{z}v|^{2}+\frac{1}{(\gamma-1)^{2}}\rho^{2-\gamma}\cdot\big(|r|^{2}+|\partial_{h}r|^{2}\big)}dxdz\Big)^{1/2}(t),
\end{align*}
and
\begin{align*}
\mathcal{E}=\mathcal{E}(t):=\Big(-\int_{\Omega}{\mathbb{S}v\cdot v}dxdz-\int_{\Omega}{\mathbb{S}(\partial_{h}v)\cdot \partial_{h}v}dxdz-\int_{\Omega}{\mathbb{S}(\partial_{z}v)\cdot\partial_{z}v}dxdz\Big)^{1/2}(t).
\end{align*}
Then it is clear that
\begin{align}\label{eq942}
\mathcal{D}\approx \|v\|_{H^{1}(\Omega)}+\|r\|_{H^{1}(\mathbb{T}^{2})},~~\mathrm{and}~~\mathcal{E}\approx \|\nabla_{h}v\|_{H^{1}(\Omega)}+\|\partial_{z}v\|_{H^{1}(\Omega)}.
\end{align}

Then, taking the $L^{2}$ inner product of \eqref{eq914} with $\rho^{L}(v-\partial_{h}^{2}v)$ and integrating by parts we get
\begin{equation}\label{eq984}
\begin{split}
&\frac{1}{2}\frac{d}{dt}\int_{\Omega}{\rho^{L}\big(|v|^{2}+|\partial_{h}v|^{2}\big)}dxdz-\int_{\Omega}{\mathbb{S}v\cdot v}dxdz-\int_{\Omega}{\mathbb{S}(\partial_{h}v)\cdot \partial_{h}v}dxdz\\
&+\int_{\Omega}{\big(u\cdot\nabla_{h}v+v\cdot\nabla_{h}u^{L}+w\cdot\partial_{z}v+(w-w^{L})\cdot\partial_{z}u^{L}\big)\cdot\rho^{L}\cdot(v-\partial_{h}^{2}v)}\\
&=\int_{\Omega}{\Big(\frac{1}{(\gamma-1)\delta}r-\Pi\Big)\cdot\mathrm{div}_{h}\,\big(\rho^{L}\cdot(v-\partial_{h}^{2}v)\big)}dxdz+\int_{\Omega}{\frac{\rho^{L}-\rho}{\rho}\cdot\mathbb{S}u\cdot(v-\partial_{h}^{2}v)}dxdz.
\end{split}
\end{equation}
From \eqref{est}, \eqref{eq943}, \eqref{eq9111}, \eqref{eq941}  and \eqref{eq942} we can verify that
\begin{equation*}
\int_{\Omega}{u\cdot\nabla_{h}v\cdot\rho^{L}\cdot(v-\partial_{h}^{2}v)}dxdz\leq C\|u\|_{H^{2}(\Omega)}\|\nabla_{h}v\|_{L^{2}(\Omega)}\|v-\partial_{h}^{2}v\|_{L^{2}(\Omega)}\leq C\mathcal{D}\big(\mathcal{D}+\mathcal{E}\big),
\end{equation*}
\begin{equation*}
\int_{\Omega}{v\cdot\nabla_{h}u^{L}\cdot\rho^{L}\cdot(v-\partial_{h}^{2}v)}dxdz\leq C\|v\|_{H^{1}(\Omega)}\|\nabla_{h}u^{L}\|_{H^{1}(\Omega)}\|v-\partial_{h}^{2}v\|_{L^{2}(\Omega)}\leq C\mathcal{D}\big(\mathcal{D}+\mathcal{E}\big),
\end{equation*}
\begin{equation*}
\int_{\Omega}{w\cdot\partial_{z}v\cdot\rho^{L}\cdot(v-\partial_{h}^{2}v)}dxdz\leq C\|w\|_{H^{2}(\Omega)}\|\partial_{z}v\|_{L^{2}(\Omega)}\|v-\partial_{h}^{2}v\|_{L^{2}(\Omega)}\leq C\mathcal{D}\big(\mathcal{D}+\mathcal{E}\big),
\end{equation*}
\begin{align*}
&\int_{\Omega}{(w-w^{L})\cdot\partial_{z}u^{L}\cdot\rho^{L}\cdot(v-\partial_{h}^{2}v)}dxdz\\
&\leq C\|w-w^{L}\|_{L^{2}(\Omega)}\|\partial_{z}u^{L}\|_{H^{2}(\Omega)}\|v-\partial_{h}^{2}v\|_{L^{2}(\Omega)}\\
&\leq C\big(\delta\|r\|_{H^{1}(\mathbb{T}^{2})}+\|v\|_{H^{1}(\Omega)}\big)\times\|\partial_{z}u^{L}\|_{H^{2}(\Omega)}\times\big(\mathcal{D}+\mathcal{E}\big)\\
&\leq C\|\partial_{z}u^{L}\|_{H^{2}(\Omega)}\mathcal{D}\big(\mathcal{D}+\mathcal{E}\big),
\end{align*}
and
\begin{equation*}
\begin{split}
\int_{\Omega}{\frac{\rho^{L}-\rho}{\rho}\cdot\mathbb{S}u\cdot(v-\partial_{h}^{2}v)}dxdz\leq C\delta\|r\|_{H^{2}(\Omega)}\|\mathbb{S}u\|_{L^{2}(\Omega)}\|v-\partial_{h}^{2}v\|_{L^{2}(\Omega)}\leq C\delta\big(\mathcal{D}+\mathcal{E}\big).
\end{split}
\end{equation*}
The estimations for the terms involved with pressure in \eqref{eq984} are more complicated. By integrating by parts, \eqref{eqL}$_{2}$  and noting that $r$, $\Pi$ are independent of $z$ we obtain
\begin{align*}
&\int_{\Omega}{\Big(\frac{1}{(\gamma-1)\delta}r-\Pi\Big)\cdot\mathrm{div}_{h}\,(\rho^{L}v)}dxdz=\int_{\Omega}{\Big(\frac{1}{(\gamma-1)\delta}r-\Pi\Big)\cdot\mathrm{div}_{h}\,(\rho^{L}u)}dxdz\\
&=-\int_{\Omega}{\Big(\frac{1}{(\gamma-1)\delta}r-\Pi\Big)\cdot\big(\partial_{t}\rho+\mathrm{div}_{h}\,\big((\rho-\rho^{L})u\big)+\partial_{z}(\rho w)\big)}dxdz\\
&=-\int_{\Omega}{\Big(\frac{1}{(\gamma-1)^{2}}r-\frac{\delta}{\gamma-1}\Pi\big)\cdot\rho^{2-\gamma}\cdot\partial_{t}r}dxdz-\int_{\Omega}{\Big(\frac{1}{(\gamma-1)\delta}r-\Pi\Big)\cdot\mathrm{div}_{h}\,\big((\rho-\rho^{L})u\big)}dxdz\\
&=-\frac{d}{dt}\int_{\Omega}{\frac{1}{2(\gamma-1)^{2}}\rho^{2-\gamma}\cdot|r|^{2}}dxdz+\frac{1}{2(\gamma-1)^{2}}\int_{\Omega}{\partial_{t}\rho^{2-\gamma}\cdot|r|^{2}}dxdz\\
&~~~+\frac{d}{dt}\int_{\Omega}{\frac{\delta}{\gamma-1}\Pi\cdot\rho^{2-\gamma}\cdot r}dxdz-\frac{\delta}{\gamma-1}\int_{\Omega}{\partial_{t}\Pi\cdot\rho^{2-\gamma}\cdot r}dxdz-\frac{\delta}{\gamma-1}\int_{\Omega}{\Pi\cdot\partial_{t}\rho^{2-\gamma}\cdot r}dxdz\\
&~~~-\int_{\Omega}{\Big(\frac{1}{(\gamma-1)\delta}r-\Pi\Big)\cdot\mathrm{div}_{h}\,\big((\rho-\rho^{L})u\big)}dxdz.
\end{align*}
Using \eqref{est}, \eqref{eq943}, \eqref{eq9112} and the estimate
\begin{align*}
\|\partial_{t}\rho\|_{H^{1}(\Omega)}\leq C\|\rho\|_{H^{2}(\Omega)}\|u\|_{H^{2}(\Omega)}+C\|\rho\|_{H^{2}(\Omega)}\|w\|_{H^{2}(\Omega)}\leq C
\end{align*}
we can verify that
\begin{align*}
&\frac{1}{2(\gamma-1)^{2}}\int_{\Omega}{\partial_{t}\rho^{2-\gamma}\cdot|r|^{2}}dxdz-\frac{\delta}{\gamma-1}\int_{\Omega}{\Pi\cdot\partial_{t}\rho^{2-\gamma}\cdot r}dxdz\\
&=\frac{2-\gamma}{2(\gamma-1)^{2}}\int_{\Omega}{\rho^{1-\gamma}\cdot\partial_{t}\rho\cdot|r|^{2}}dxdz-\frac{\delta(2-\gamma)}{\gamma-1}\int_{\Omega}{\Pi\cdot\rho^{1-\gamma}\cdot\partial_{t}\rho\cdot r}dxdz\\
&\leq C\big(\|r\|^{2}_{H^{1}(\Omega)}+\delta\|\Pi\|_{H^{1}(\Omega)}\|r\|_{H^{1}(\Omega)}\big)\|\partial_{t}\rho\|_{H^{1}(\Omega)}\leq C\big(\mathcal{D}^{2}+\delta\mathcal{D}\big),
\end{align*}
\begin{align*}
-\frac{\delta}{\gamma-1}\int_{\Omega}{\partial_{t}\Pi\cdot\rho^{2-\gamma}\cdot r}dxdz\leq C\delta\|\partial_{t}\Pi\|_{L^{2}(\Omega)}\|r\|_{L^{2}(\mathbb{T}^{2})}\leq C\delta\|\partial_{t}\Pi\|_{L^{2}(\mathbb{T}^{2})}\mathcal{D},
\end{align*}
and
\begin{align*}
&-\int_{\Omega}{\Big(\frac{1}{(\gamma-1)\delta}r-\Pi\Big)\cdot\mathrm{div}_{h}\,\big((\rho-\rho^{L})u\big)}dxdz\\
&\leq C\big(\frac{1}{\delta}\|r\|_{L^{2}(\Omega)}+\|\Pi\|_{L^{2}(\Omega)}\big)\times\big(\|\rho-\rho^{L}\|_{L^{4}(\Omega)}\|\mathrm{div}_{h}\,u\|_{L^{4}(\Omega)}+\|\nabla_{h}\rho\|_{L^{2}(\Omega)}\|u\|_{L^{\infty}(\Omega)}\big)\\
&\leq C\big(\frac{1}{\delta}\|r\|_{L^{2}(\mathbb{T}^{2})}+\|\Pi\|_{L^{2}(\mathbb{T}^{2})}\big)\times\big(\delta\|r\|_{L^{4}(\mathbb{T}^{2})}+\delta\|\partial_{h}r\|_{L^{2}(\mathbb{T}^{2})}\big)\leq C\big(\mathcal{D}^{2}+\delta\mathcal{D}\big).
\end{align*}
Thus, combining the estimates above we have
\begin{equation}\label{eq944}
\begin{split}
&\int_{\Omega}{\Big(\frac{1}{(\gamma-1)\delta}r-\Pi\Big)\cdot\mathrm{div}_{h}\,(\rho^{L}v)}dxdz\\
&\leq-\frac{d}{dt}\int_{\Omega}{\frac{1}{2(\gamma-1)^{2}}\rho^{2-\gamma}\cdot|r|^{2}}dxdz+\frac{d}{dt}\int_{\Omega}{\frac{\delta}{\gamma-1}\Pi\cdot\rho^{2-\gamma}\cdot r}dxdz\\
&~~~~+C\big(\mathcal{D}^{2}+\delta\mathcal{D}\big)+C\delta\|\partial_{t}\Pi\|_{L^{2}(\mathbb{T}^{2})}\mathcal{D}.
\end{split}
\end{equation}
In similar fashion to the previous discussions, we have
\begin{equation}\label{eq945}
\begin{split}
&-\int_{\Omega}{\Big(\frac{1}{(\gamma-1)\delta}r-\Pi\Big)\cdot\mathrm{div}_{h}\,(\rho^{L}\partial_{h}^{2}v)}dxdz=\int_{\Omega}{\partial_{h}\Big(\frac{1}{(\gamma-1)\delta}r-\Pi\Big)\cdot\partial_{h}\mathrm{div}_{h}\,(\rho^{L}u)}dxdz\\
&=-\int_{\Omega}{\partial_{h}\Big(\frac{1}{(\gamma-1)\delta}r-\Pi\Big)\cdot\big(\partial_{t}\partial_{h}\rho+\partial_{h}\mathrm{div}_{h}\,\big((\rho-\rho^{L})u\big)+\partial_{h}\partial_{z}(\rho w)\big)}dxdz\\
&=-\int_{\Omega}{\partial_{h}\Big(\frac{1}{(\gamma-1)\delta}r-\Pi\Big)\cdot\Big(\frac{\delta}{\gamma-1}\cdot\partial_{t}\big(\rho^{2-\gamma}\partial_{h}r\big)+\partial_{h}\mathrm{div}_{h}\,\big((\rho-\rho^{L})u\big)\Big)}dxdz,
\end{split}
\end{equation}
where, by using the fact
\begin{equation*}
\begin{split}
&\|\partial_{t}\rho\|_{L^{\infty}(\Omega)}\leq C\big(\|\mathrm{div}_{h}(\rho u)\|_{L^{\infty}(\Omega)}+\|\partial_{z}(\rho w)\|_{L^{\infty}(\Omega)}\big)\\ &\leq C\big(\|u\|_{H^{2}(\Omega)}\|\nabla_{h}\rho\|_{L^{\infty}(\Omega)}+\|\mathrm{div}_{h}\,u\|_{H^{2}(\Omega)}+\|\partial_{z}\rho\|_{L^{\infty}(\Omega)}\|w\|_{H^{2}(\Omega)}+\|\partial_{z}w\|_{L^{\infty}(\Omega)}\big)\\
&\leq C\Big(\mathcal{A}^{2}+\mathcal{B}+\mathcal{A}+\mathcal{A}^{4}+\|(\xi^{-1},\partial_{z}(\xi^{-1}))\|_{L^{\infty}(\Omega)}\times\|(\mathrm{div}_{h}(\xi u),\delta u\cdot \nabla_{h}r,w)\|_{L^{\infty}(\Omega)}\big)\Big)\\
&\leq C\big(\mathcal{B}+1\big)
\end{split}
\end{equation*}
we have
\begin{align}\nonumber
&-\int_{\Omega}{\partial_{h}\Big(\frac{1}{(\gamma-1)\delta}r-\Pi\Big)\cdot\frac{\delta}{\gamma-1}\cdot\partial_{t}\big(\rho^{2-\gamma}\partial_{h}r\big)}dxdz\\\nonumber
&=-\frac{d}{dt}\int_{\Omega}{\frac{1}{2(\gamma-1)^{2}}\rho^{2-\gamma}\cdot|\partial_{h}r|^{2}}dxdz+\frac{d}{dt}\int_{\Omega}{\frac{\delta}{\gamma-1}\cdot\partial_{h}\Pi\cdot\rho^{2-\gamma}\cdot\partial_{h} r}dxdz\\\nonumber
&~~~~-\frac{1}{2(\gamma-1)^{2}}\int_{\Omega}{\partial_{t}\rho^{2-\gamma}\cdot|\partial_{h}r|^{2}}dxdz-\frac{\delta}{\gamma-1}\int_{\Omega}{\partial_{t}\partial_{h}\Pi\cdot\rho^{2-\gamma}\cdot \partial_{h}r}dxdz\\\label{eq985}
&\leq -\frac{d}{dt}\int_{\Omega}{\frac{1}{2(\gamma-1)^{2}}\rho^{2-\gamma}\cdot|\partial_{h}r|^{2}}dxdz+\frac{d}{dt}\int_{\Omega}{\frac{\delta}{\gamma-1}\cdot\partial_{h}\Pi\cdot\rho^{2-\gamma}\cdot\partial_{h} r}dxdz\\\nonumber
&~~~~+C\|\partial_{t}\rho\|_{L^{\infty}(\Omega)}\|\partial_{h}r\|^{2}_{L^{2}(\Omega)}+C\delta\|\partial_{t}\partial_{h}\Pi\|_{L^{2}(\Omega)}\|\partial_{h}r\|_{L^{2}(\Omega)}\\\nonumber
&\leq -\frac{d}{dt}\int_{\Omega}{\frac{1}{2(\gamma-1)^{2}}\rho^{2-\gamma}\cdot|\partial_{h}r|^{2}}dxdz+\frac{d}{dt}\int_{\Omega}{\frac{\delta}{\gamma-1}\cdot\partial_{h}\Pi\cdot\rho^{2-\gamma}\cdot\partial_{h} r}dxdz\\\nonumber
&~~~~+C\big(\mathcal{B}+1\big)\mathcal{D}^{2}+C\delta\|\partial_{t}\partial_{h}\Pi\|_{L^{2}(\mathbb{T}^{2})}\mathcal{D}.\nonumber
\end{align}
To estimate the another term in \eqref{eq945}, we split
\begin{equation}\label{eq946}
\begin{split}
&-\int_{\Omega}{\partial_{h}\Big(\frac{1}{(\gamma-1)\delta}r-\Pi\Big)\cdot\partial_{h}\mathrm{div}_{h}\,\big((\rho-\rho^{L})u\big)}dxdz\\
&=-\int_{\Omega}{\partial_{h}\Big(\frac{1}{(\gamma-1)\delta}r-\Pi\Big)\cdot\big((\rho-\rho^{L})\mathrm{div}_{h}\,\partial_{h}u+\partial_{h}\rho\cdot\mathrm{div}_{h}\,u+\partial_{h}u\cdot\nabla_{h}\rho\big)}dxdz\\
&~~~-\frac{1}{(\gamma-1)\delta}\int_{\Omega}{\partial_{h}r\cdot u\cdot\nabla_{h}\partial_{h}\rho}dxdz+\int_{\Omega}{\partial_{h}\Pi\cdot u\cdot\nabla_{h}\partial_{h}\rho}dxdz,
\end{split}
\end{equation}
where, from \eqref{eq943} and \eqref{eq9112} we can check that
\begin{equation*}
\begin{split}
&-\int_{\Omega}{\partial_{h}\Big(\frac{1}{(\gamma-1)\delta}r-\Pi\Big)\cdot\big((\rho-\rho^{L})\mathrm{div}_{h}\,\partial_{h}u+\partial_{h}\rho\cdot\mathrm{div}_{h}\,u+\partial_{h}u\cdot\nabla_{h}\rho\big)}dxdz\\
&\leq C\big(\frac{1}{\delta}\|\partial_{h}r\|_{L^{2}(\Omega)}+\|\partial_{h}\Pi\|_{L^{2}(\Omega)}\big)\times\big(\delta\|r\|_{L^{4}(\Omega)}\|\partial_{h}^{2}u\|_{L^{4}(\Omega)}+\|\partial_{h}\rho\|_{L^{2}(\Omega)}\|\partial_{h}u\|_{L^{\infty}(\Omega)}\big)\\
&\leq C\big(\frac{1}{\delta}\mathcal{D}+1\big)\times\delta\mathcal{D}\mathcal{B}\leq C\big(\mathcal{B}\mathcal{D}^{2}+\delta\mathcal{B}\mathcal{D}\big).
\end{split}
\end{equation*}
Furthermore, by using integration by parts, \eqref{est} and \eqref{eq9112} we have
\begin{align*}
&-\frac{1}{(\gamma-1)\delta}\int_{\Omega}{\partial_{h}r\cdot u\cdot\nabla_{h}\partial_{h}\rho}dxdz\\
&=-\frac{1}{(\gamma-1)^{2}}\int_{\Omega}{\rho^{2-\gamma}\cdot u\cdot\nabla_{h}\partial_{h}r\cdot\partial_{h}r}dxdz-\frac{\delta(2-\gamma)}{(\gamma-1)^{3}}\int_{\Omega}{\partial_{h}r\cdot u\cdot\rho^{3-2\gamma}\nabla_{h}r\cdot\partial_{h}r}dxdz\\
&=\frac{1}{(\gamma-1)^{2}}\int_{\Omega}{\nabla_{h}\rho^{2-\gamma}\cdot u\cdot\frac{1}{2}|\partial_{h}r|^{2}}dxdz+\frac{1}{(\gamma-1)^{2}}\int_{\Omega}{\rho^{2-\gamma}\cdot \mathrm{div}_{h}\,u\cdot\frac{1}{2}|\partial_{h}r|^{2}}dxdz\\
&~~~~-\frac{\delta(2-\gamma)}{(\gamma-1)^{3}}\int_{\Omega}{\partial_{h}r\cdot u\cdot\rho^{3-2\gamma}\nabla_{h}r\cdot\partial_{h}r}dxdz\\
&\leq C\big(\delta\|\nabla_{h}r\|_{L^{\infty}(\Omega)}\|u\|_{L^{\infty}(\Omega)}+\|\mathrm{div}_{h}\,u\|_{L^{\infty}(\Omega)}\big)\|\partial_{h}r\|^{2}_{L^{2}(\Omega)}\\
&\leq C\big(\delta\mathcal{A}^{2}+\mathcal{B}\big)\mathcal{D}^{2}\leq C\big(\mathcal{D}^{2}+\mathcal{B}\mathcal{D}^{2}\big),
\end{align*}
and
\begin{equation*}
\begin{split}
&\int_{\Omega}{\partial_{h}\Pi\cdot u\cdot\nabla_{h}\partial_{h}\rho}dxdz=-\int_{\Omega}{\nabla_{h}\partial_{h}\Pi\cdot u\cdot\partial_{h}\rho}dxdz-\int_{\Omega}{\partial_{h}\Pi\cdot \mathrm{div}_{h}\,u\cdot\partial_{h}\rho}dxdz\\
&\leq C\big(\|\nabla_{h}\partial_{h}\Pi\|_{L^{2}(\Omega)}\|u\|_{L^{\infty}(\Omega)}+\|\partial_{h}\Pi\|_{L^{4}(\Omega)}\|\partial_{h}u\|_{L^{4}(\Omega)}\big)\times\delta\|\partial_{h}r\|_{L^{2}(\mathbb{T}^{2})}\leq C\delta \mathcal{D}.
\end{split}
\end{equation*}
Plugging the above estimates into \eqref{eq946} we get
\begin{equation}\label{eq986}
-\int_{\Omega}{\partial_{h}\Big(\frac{1}{(\gamma-1)\delta}r-\Pi\Big)\cdot\partial_{h}\mathrm{div}_{h}\,\big((\rho-\rho^{L})u\big)}dxdz\leq C\big(\mathcal{B}\mathcal{D}^{2}+\delta\mathcal{B}\mathcal{D}+\mathcal{D}^{2}+\delta\mathcal{D}\big).
\end{equation}
As a result, combining \eqref{eq945}, \eqref{eq985} and \eqref{eq986} we obtain
\begin{align*}
&-\int_{\Omega}{\Big(\frac{1}{(\gamma-1)\delta}r-\Pi\Big)\cdot\mathrm{div}_{h}\,(\rho^{L}\partial_{h}^{2}v)}dxdz\\
&\leq -\frac{d}{dt}\int_{\Omega}{\frac{1}{2(\gamma-1)^{2}}\rho^{2-\gamma}\cdot|\partial_{h}r|^{2}}dxdz+\frac{d}{dt}\int_{\Omega}{\frac{\delta}{\gamma-1}\cdot\partial_{h}\Pi\cdot\rho^{2-\gamma}\cdot\partial_{h} r}dxdz\\
&~~~~+C\big(\mathcal{B}\mathcal{D}^{2}+\mathcal{D}^{2} +\delta\|\partial_{t}\partial_{h}\Pi\|_{L^{2}(\mathbb{T}^{2})}\mathcal{D}+\delta\mathcal{B}\mathcal{D}+\delta\mathcal{D}\big).
\end{align*}
Finally, collecting the discussions above we conclude with
\begin{equation}\label{eq9413}
\begin{split}
&\frac{1}{2}\frac{d}{dt}\int_{\Omega}{\rho^{L}\big(|v|^{2}+|\partial_{h}v|^{2}\big)+\frac{1}{(\gamma-1)^{2}}\rho^{2-\gamma}\cdot\big(|r|^{2}+|\partial_{h}r|^{2}\big)}dxdz\\
&~~~-\int_{\Omega}{\mathbb{S}v\cdot v}dxdz-\int_{\Omega}{\mathbb{S}(\partial_{h}v)\cdot \partial_{h}v}dxdz\\
&\leq \frac{d}{dt}\int_{\Omega}{\frac{\delta}{\gamma-1}\rho^{2-\gamma}\cdot\big(\Pi\cdot r+\partial_{h}\Pi\cdot\partial_{h}r\big)}dxdz+C\|\partial_{z}u^{L}\|_{H^{2}(\Omega)}\big(\mathcal{D}^{2}+\mathcal{D}\mathcal{E}\big)\\
&~~~+C\big(\mathcal{D}^{2}+\mathcal{D}\mathcal{E}+\mathcal{B}\mathcal{D}^{2}+\delta\mathcal{D}+\delta\mathcal{E}+\delta\mathcal{B}\mathcal{D}+\delta\|\partial_{t}\Pi\|_{H^{1}(\mathbb{T}^{2})}\mathcal{D}\big).
\end{split}
\end{equation}

To derive the estimates for $\partial_{z}v$, we take the $L^{2}$ inner product of \eqref{eq914} with $-\partial_{z}^{2}v$ and integrate by parts to obtain
\begin{align*}
&\frac{1}{2}\frac{d}{dt}\int_{\Omega}|\partial_{z}v|^{2}dxdz-\int_{\Omega}{\mathbb{S}(\partial_{z}v)\cdot\partial_{z}v}dxdz\\
&=\int_{\Omega}{\big(u\cdot\nabla_{h}v+v\cdot\nabla_{h}u^{L}+w\cdot\partial_{z}v+(w-w^{L})\cdot\partial_{z}u^{L}\big)\cdot\partial_{z}^{2}v}dxdz\\
&~~~~+\int_{\Omega}{\Big(1-\frac{1}{\rho^{L}}\Big)\mathbb{S}v\cdot\partial_{z}^{2}v}dxdz+\int_{\Omega}{\Big(\frac{1}{\rho^{L}}-\frac{1}{\rho}\Big)\mathbb{S}u\cdot\partial_{z}^{2}v}dxdz.
\end{align*}
By \eqref{est}, \eqref{eq8232}, \eqref{eq943}, \eqref{eq9112} and \eqref{eq941} we have
\begin{align*}
&\int_{\Omega}{\big(u\cdot\nabla_{h}v+v\cdot\nabla_{h}u^{L}+w\cdot\partial_{z}v+(w-w^{L})\cdot\partial_{z}u^{L}\big)\cdot\partial_{z}^{2}v}dxdz\\
&\leq C\Big(\|u\|_{H^{2}(\Omega)}\|\nabla_{h}v\|_{L^{2}(\Omega)}+\|v\|_{H^{1}(\Omega)}\|\nabla_{h}u^{L}\|_{H^{1}(\Omega)}+\|w\|_{H^{2}(\Omega)}\|\partial_{z}v\|_{L^{2}(\Omega)}\\
&~~~~+\|w-w^{L}\|_{L^{2}(\Omega)}\|\partial_{z}u^{L}\|_{H^{2}(\Omega)}\Big)\|\partial_{z}^{2}v\|_{L^{2}(\Omega)}\\
&\leq C\big(\mathcal{A}\mathcal{D}+\mathcal{D}+\big(\mathcal{A}+\mathcal{A}^{4}\big)\mathcal{D}+\|\partial_{z}u^{L}\|_{H^{2}(\Omega)}\mathcal{D}\big)\mathcal{E}\\
&\leq C\big(\mathcal{D}\mathcal{E}+\|\partial_{z}u^{L}\|_{H^{2}(\Omega)}\mathcal{D}\mathcal{E}\big),
\end{align*}
and
\begin{align*}
\int_{\Omega}{\Big(\frac{1}{\rho^{L}}-\frac{1}{\rho}\Big)\mathbb{S}u\cdot\partial_{z}^{2}v}dxdz\leq C\delta\|r\|_{L^{4}(\Omega)}\|\mathbb{S}u\|_{L^{4}(\Omega)}\|\partial_{z}^{2}v\|_{L^{2}(\Omega)}\leq C\mathcal{D}\mathcal{B}\mathcal{E}.
\end{align*}
Then, noting that $1-\frac{1}{\rho^{L}}<0$, we have
\begin{align*}
&\int_{\Omega}{\Big(1-\frac{1}{\rho^{L}}\Big)\mathbb{S}v\cdot\partial_{z}^{2}v}dxdz\leq\int_{\Omega}{\Big(1-\frac{1}{\rho^{L}}\Big)\cdot\big(\mu\triangle_{h}v+\nu\nabla_{h}\mathrm{div}_{h}\,v\big)\cdot\partial_{z}^{2}v}dxdz\\
&=\mu\int_{\Omega}{\Big(1-\frac{1}{\rho^{L}}\Big)\cdot|\nabla_{h}\partial_{z}v|^{2}}dxdz+\mu\int_{\Omega}{\partial_{z}\Big(\frac{1}{\rho^{L}}\Big)\cdot\triangle_{h}v\cdot\partial_{z}v}dxdz\\
&~~~~+\nu\int_{\Omega}{\Big(1-\frac{1}{\rho^{L}}\Big)\cdot|\mathrm{div}_{h}\,\partial_{z}v|^{2}}dxdz+\nu\int_{\Omega}{\partial_{z}\Big(\frac{1}{\rho^{L}}\Big)\cdot\nabla_{h}\mathrm{div}_{h}\,v\cdot\partial_{z}v}dxdz\\
&\leq\mu\int_{\Omega}{\partial_{z}\Big(\frac{1}{\rho^{L}}\Big)\cdot\triangle_{h}v\cdot\partial_{z}v}dxdz+\nu\int_{\Omega}{\partial_{z}\Big(\frac{1}{\rho^{L}}\Big)\cdot\nabla_{h}\mathrm{div}_{h}\,v\cdot\partial_{z}v}dxdz\\
&\leq C\|\partial_{h}^{2}v\|_{L^{2}(\Omega)}\|\partial_{z}v\|_{L^{2}(\Omega)}\leq C\mathcal{E}\mathcal{D}.
\end{align*}
Collecting the estimates above we get
\begin{equation}\label{eq9412}
\frac{1}{2}\frac{d}{dt}\int_{\Omega}{|\partial_{z}v|^{2}}dxdz-\int_{\Omega}{\mathbb{S}(\partial_{z}v)\cdot\partial_{z}v}dxdz\leq C\big(\mathcal{D}\mathcal{E}+\|\partial_{z}u^{L}\|_{H^{2}(\Omega)}\mathcal{D}\mathcal{E}+\mathcal{B}\mathcal{D}\mathcal{E}\big).
\end{equation}

Now, combining \eqref{eq9413} and \eqref{eq9412} and using Young's inequality we arrive at
\begin{equation*}
\begin{split}
&\frac{1}{2}\frac{d}{dt}\mathcal{D}^{2}+\mathcal{E}^{2}\leq\frac{d}{dt}\int_{\Omega}{\frac{\delta}{\gamma-1}\rho^{2-\gamma}\cdot\big(\Pi\cdot r+\partial_{h}\Pi\cdot\partial_{h}r\big)}dxdz+C\|\partial_{z}u^{L}\|_{H^{2}(\Omega)}\big(\mathcal{D}^{2}+\mathcal{D}\mathcal{E}\big)\\
&~~~~+C\big(\mathcal{D}^{2}+\mathcal{D}\mathcal{E}+\mathcal{B}\mathcal{D}^{2}+\delta\mathcal{D}+\delta\mathcal{E}+\delta\mathcal{B}\mathcal{D}+\delta\|\partial_{t}\Pi\|_{H^{1}(\mathbb{T}^{2})}\mathcal{D}+\mathcal{B}\mathcal{D}\mathcal{E}\big)\\
&\leq\frac{d}{dt}\int_{\Omega}{\frac{\delta}{\gamma-1}\rho^{2-\gamma}\cdot\big(\Pi\cdot r+\partial_{h}\Pi\cdot\partial_{h}r\big)}dxdz+C\Big(\|\partial_{z}u^{L}\|_{H^{2}(\Omega)}+\|\partial_{z}u^{L}\|_{H^{2}(\Omega)}^{2}+1\\
&~~~~+\mathcal{B}+\mathcal{B}^{2}\Big)\mathcal{D}^{2}+C\delta^{2}\big(1+\mathcal{B}^{2}+\|\partial_{t}\Pi\|_{H^{1}(\mathbb{T}^{2})}^{2}\big)+\frac{1}{2}\mathcal{E}^{2}.
\end{split}
\end{equation*}
Thus, integrating the above inequality over $(0,t)$ and noting that
\begin{equation*}
\begin{split}
&\int^{t}_{0}{\frac{d}{dt}\int_{\Omega}{\frac{\delta}{\gamma-1}\rho^{2-\gamma}\cdot\big(\Pi\cdot r+\partial_{h}\Pi\cdot\partial_{h}r\big)(\tau)}dxdz}d\tau\\
&=\frac{\delta}{\gamma-1}\int_{\Omega}{\rho^{2-\gamma}\cdot\big(\Pi\cdot r+\partial_{h}\Pi\cdot\partial_{h}r\big)(t)}dxdz-\frac{\delta}{\gamma-1}\int_{\Omega}{\rho^{2-\gamma}\cdot\big(\Pi\cdot r+\partial_{h}\Pi\cdot\partial_{h}r\big)|_{t=0}}dxdz\\
&\leq C\delta\|\Pi(t)\|_{H^{1}(\mathbb{T}^{2})}\|r(t)\|_{H^{1}(\mathbb{T}^{2})}+C\delta\|\Pi|_{t=0}\|_{H^{1}(\mathbb{T}^{2})}\|r_{0}\|_{H^{1}(\mathbb{T}^{2})}\\
&\leq C\delta\mathcal{D}(t)+C\delta\mathcal{D}(0)\leq \frac{1}{4}\mathcal{D}^{2}(t)+\frac{1}{2}\mathcal{D}^{2}(0)+C\delta^{2}
\end{split}
\end{equation*}
in view of \eqref{eq9112}, we get
\begin{align*}
&\frac{1}{2}\mathcal{D}^{2}(t)+\frac{1}{2}\int^{t}_{0}{\mathcal{E}^{2}}(s)ds\\
&\leq \frac{1}{2}\mathcal{D}^{2}(0)+\frac{1}{4}\mathcal{D}^{2}(t)+\frac{1}{2}\mathcal{D}^{2}(0)+C\delta^{2}
+C\delta^{2}\int^{t}_{0}{\big(1+\mathcal{B}^{2}+\|\partial_{t}\Pi\|_{H^{1}(\Omega)}^{2}\big)}ds\\
&~~~+C\int^{t}_{0}{\big(\|\partial_{z}u^{L}\|_{H^{2}(\Omega)}+\|\partial_{z}u^{L}\|_{H^{2}(\Omega)}^{2}+1+\mathcal{B}+\mathcal{B}^{2}\big)\mathcal{D}^{2}(s)}ds\\
&\leq \mathcal{D}^{2}(0)+\frac{1}{4}\mathcal{D}^{2}(t)+C\delta^{2}+C\int^{t}_{0}{\big(\|\partial_{z}u^{L}\|_{H^{2}(\Omega)}+\|\partial_{z}u^{L}\|_{H^{2}(\Omega)}^{2}+1+\mathcal{B}+\mathcal{B}^{2}\big)\mathcal{D}^{2}(s)}ds.
\end{align*}
Thus, by applying Gronwall's inequality we arrive at
\begin{equation*}
\begin{split}
\mathcal{D}^{2}(t)+2\int^{t}_{0}{\mathcal{E}^{2}}(s)ds\leq e^{C\int^{t}_{0}{(\|\partial_{z}u^{L}\|_{H^{2}(\Omega)}+\|\partial_{z}u^{L}\|_{H^{2}(\Omega)}^{2}+1+\mathcal{B}+\mathcal{B}^{2})(s)}ds}\big(4\mathcal{D}^{2}(0)+C\delta^{2}\big),~~\forall\,t\in[0,T],
\end{split}
\end{equation*}
where the constants $C$ above are independent of $\delta$. Finally, by using \eqref{est}, \eqref{eq983}, \eqref{eq9111} and \eqref{eq942} we can obtain the estimate \eqref{eqWell}. Now we have completed the proof of Theorem \ref{th2}.

\section{Low Mach number limit with ill-prepared data}\label{sec4}
In this section we study briefly the low Mach number singular limit of the compressible primitive equations with ill-prepared data. Different from the discussions in the previous section, the ill-prepared data allow the propagation of the fast-oscillating acoustic waves in the solutions.
These fast waves shall not decay in periodic area. As a result, the solutions of \eqref{eq2} are expected to converge weakly as the Mach number tends to zero.

Let $\{(r,u)\}_{\delta}$ be the solutions of \eqref{eq2} given by Theorem \ref{th1}, and let $\rho$ be defined by \eqref{eqrho}. Then $\{(\rho,u)\}_{\delta}$ are strong solutions of \eqref{e1}. By \eqref{est}, \eqref{eq-est-w} and the differential mean value theorem we have
\begin{equation}\label{eq962}
\|\rho-\rho^{L}\|_{L^{\infty}(0,T;H^{3}(\Omega))}\leq C\delta,
\end{equation}
and as $\delta$ goes to zero, up to a subsequence,
\begin{equation}\label{eq961}
(u,w)\to(u^{L},w^{L})~~\mathrm{wakely*~in}~~L^{\infty}(0,T;H^{2}(\Omega))
\end{equation}
for some $(u^{L},w^{L})\in L^{\infty}(0,T;H^{2}(\Omega))$. Thus, by taking $\delta\to0$ (or up to a subsequence) in the weak formulation of the density equation \eqref{e1}$_{1}$ we get
\begin{equation*}
\mathrm{div}_{h}(\rho^{L}u^{L})+\partial_{z}(\rho^{L}w^{L})=0.
\end{equation*}
Furthermore, after applying $\partial_{z}$ and $\mathrm{curl}_{h}$ to \eqref{eq2}$_{2}$ to annihilate the singular pressure term, we can verify that
\begin{equation*}
\|(\partial_{t}\partial_{z}u,\partial_{t}\mathrm{curl}_{h}u)\|_{L^{2}(0,T;L^{2}(\Omega))}\leq C.
\end{equation*}
Thus, by \eqref{est}, Lions-Aubin's lemma and \eqref{eq961} we get
\begin{equation}\label{eq963}
(\partial_{z}u,\mathrm{curl}_{h}u)\to(\partial_{z}u^{L},\mathrm{curl}_{h}u^{L})~~\mathrm{strongly~in}~~L^{\infty}(0,T;L^{2}(\Omega)),
\end{equation}
as $\delta\to0$, up to a subsequence.

Now we show that $(\rho^{L},u^{L},w^{L})$ solves the equation \eqref{eqL}$_{1}$. As stated before, as $\delta\to0$, the sequence $\{u\}_{\delta}$ shall converge in weak sense. Thus, it's more convenient to identify the limiting system in the integral formula of \eqref{eq2}$_{2}$. The convergences in \eqref{eq962}, \eqref{eq961} and \eqref{eq963} are enough to achieve this. In the rest of this section, without further mentioned, all the convergences hold in the sense of
a subsequence of $\delta$, as $\delta\to0$.

Letting $\phi\in C_{c}^{\infty}(\overline{\Omega}\times[0,T);\mathbb{R}^{2})$ such that $\mathrm{div}_{h}\,\phi=0$, taking the inner product of \eqref{eq2}$_{2}$ with $\phi$ and integrating by parts we have
\begin{equation}\label{eq965}
\begin{split}
&-\int_{0}^{T}\int_{\Omega}{u\cdot\partial_{t}\phi}dxdzdt+\int_{0}^{T}\int_{\Omega}{u\cdot\nabla_{h}u\cdot\phi}dxdzdt+\int_{0}^{T}\int_{\Omega}{w\cdot\partial_{z}u\cdot\phi}dxdzdt\\
&+\int_{0}^{T}\int_{\Omega}{\frac{1}{Ro}u^{\bot}\cdot\phi}dxdzdt=\int_{0}^{T}\int_{\Omega}{\frac{1}{\rho}\mathbb{S}u\cdot\phi}dxdzdt+\int_{\Omega}{u_{0}\cdot\phi|_{t=0}}dxdz.
\end{split}
\end{equation}
From \eqref{eq9251}, \eqref{eq962}, \eqref{eq961} and \eqref{eq963} we get, after passing to the limit $\delta\to0$ in \eqref{eq965}, that
\begin{equation}\label{eq966}
\begin{split}
&-\int_{0}^{T}\int_{\Omega}{u^{L}\cdot\partial_{t}\phi}dxdzdt+\lim\limits_{\delta\to0}\int_{0}^{T}\int_{\Omega}{u\cdot\nabla_{h}u\cdot\phi}dxdzdt+\int_{0}^{T}\int_{\Omega}{w^{L}\cdot\partial_{z}u^{L}\cdot\phi}dxdzdt\\
&+\int_{0}^{T}\int_{\Omega}{\frac{1}{Ro}(u^{L})^{\bot}\cdot\phi}dxdzdt=\int_{0}^{T}\int_{\Omega}{\frac{1}{\rho^{L}}\mathbb{S}u^{L}\cdot\phi}dxdzdt+\int_{\Omega}{u_{0}^{L}\cdot\phi|_{t=0}}dxdz.
\end{split}
\end{equation}
To identify the limit of the integral involved with $u\cdot\nabla_{h}u$, we split
\begin{align}\label{eq9281}
\int_{0}^{T}\int_{\Omega}{u\cdot\nabla_{h}u\cdot\phi}dxdzdt=\int_{0}^{T}\int_{\Omega}{\frac{1}{2}\nabla_{h}|u|^{2}\cdot\phi}dxdzdt+\int_{0}^{T}\int_{\Omega}{\mathrm{curl}_{h}u\cdot u^{\bot}\cdot\phi}dxdzdt,
\end{align}
where we have used the identity
\begin{equation*}
u\cdot\nabla_{h}u=\frac{1}{2}\nabla_{h}|u|^{2}+\mathrm{curl}_{h}u\cdot u^{\bot}.
\end{equation*}
Noting that the first term in the right side of \eqref{eq9281} disappears after using integration by parts, we get, by virtue of \eqref{eq961} and \eqref{eq963},
\begin{equation}\label{eq967}
\begin{split}
&\int_{0}^{T}\int_{\Omega}{u\cdot\nabla_{h}u\cdot\phi}dxdzdt\to\int_{0}^{T}\int_{\Omega}{\mathrm{curl}_{h}u^{L}\cdot (u^{L})^{\bot}\cdot\phi}dxdzdt\\
&=\int_{0}^{T}\int_{\Omega}{u^{L}\cdot\nabla_{h}u^{L}\cdot\phi}dxdzdt-\int_{0}^{T}\int_{\Omega}{\frac{1}{2}\nabla_{h}|u^{L}|^{2}\cdot\phi}dxdzdt\\
&=\int_{0}^{T}\int_{\Omega}{u^{L}\cdot\nabla_{h}u^{L}\cdot\phi}dxdzdt,~~~~\mathrm{as}~~\delta\to0.
\end{split}
\end{equation}
Plugging \eqref{eq967} into \eqref{eq966} and taking notice of
\begin{equation*}
\int_{\Omega}{u_{0}^{L}\cdot\phi|_{t=0}}dxdz=\int_{\Omega}{\mathbf{P}u_{0}^{L}\cdot\phi|_{t=0}}dxdz,
\end{equation*}
we conclude that $(\rho^{L},u^{L},w^{L})$ solves \eqref{eqL}$_{1}$ with $u^{L}|_{t=0}=\mathbf{P}u^{L}_{0}$ for some function $\Pi$. Moreover, as
\begin{equation*}
\partial_{z}\triangle^{-1}_{h}\mathrm{div}_{h}\big(\partial_{t}u+u\cdot\nabla_{h} u+w\cdot\partial_{z}u+\frac{1}{Ro}u^{\bot}-\frac{1}{\rho}\mathbb{S}u\big)=-\frac{1}{(\gamma-1)\delta}\partial_{z}r=0,
\end{equation*}
the weak limit $\Pi$ of that term is also independent of $z$. Finally, since the strong solution of \eqref{eqL} is unique, the convergences above occur without restricting to a subsequence. Now we have completed the proof of Theorem \ref{th3}.

\begin{appendix}
\section{Appendix}
We state here how we derive the non-dimensional equations \eqref{e1} from the original compressible primitive system \eqref{CPE}. We mention that the procedure of nondimensionalization below is similar to that in \cite{B1} (Chapter 6.1).

First, from \eqref{CPE}$_{3}$ we can verify that
\begin{equation*}
\rho(x,z,t)=\bigg(\rho_{h}(x,t)-\frac{(\gamma-1)g}{\gamma}z\bigg)^{\frac{1}{\gamma-1}}
\end{equation*}
for some function $\rho_{h}(x,t)$. We assume that $\rho_{h}(x,t)$ is given by a mean $\rho_{h,m}$ and the deviation $r$:
\begin{equation*}
\rho(x,t)=\rho_{h,m}+r(x,t),
\end{equation*}
and let $\rho_{b}$ be the background density. Then we have
\begin{equation*}
\rho(x,z,t)=\rho_{b}\bigg(\frac{\rho_{h,m}}{\rho^{\gamma-1}_{b}}+\frac{r(x,t)}{\rho_{b}^{\gamma-1}}-\frac{(\gamma-1)g}{\gamma\rho_{b}^{\gamma-1}}z\bigg)^{\frac{1}{\gamma-1}}.
\end{equation*}
Without loss of generality, we assume that $\rho_{h,m}=\rho_{b}^{\gamma-1}$. So that the constant $\rho_{h,m}/\rho^{\gamma-1}_{b}$ in the above expression equals to $1$. We now introduce the following scales for length, time, velocity and density:
\begin{equation*}
\begin{split}
&L:~\mathrm{typical~length~scale};~~U:~\mathrm{typical~velocity~scale};\\
&T_{e}=\frac{L}{U}:~\mathrm{eddy~turnover~time};~~\rho_{b}:\mathrm{typical~density~scale};\\
&r_{p}: \mathrm{mean~size~of~the~density~perturbation}.
\end{split}
\end{equation*}
With this choice of scales we introduce the following nondimensional variables:
\begin{equation}\label{eq9143}
x'=\frac{x}{L},~~z'=\frac{z}{L},~~t'=\frac{t}{T_{e}},~~u'=\frac{u}{U},~~w'=\frac{w}{U},~~\rho'=\frac{\rho}{\rho_{b}},~~r'=\frac{r}{r_{p}},
\end{equation}
and the following parameters:
\begin{equation}\label{eq9142}
\begin{split}
&Ro=\frac{U}{fL},~~\delta=\frac{U}{\sqrt{\gamma\rho_{b}^{\gamma-1}}},~~\theta=\frac{(\gamma-1)gL}{\gamma\rho^{\gamma-1}_{b}},~~\mu'=\frac{\mu}{\rho_{b}LU},~~\nu'=\frac{\nu}{\rho_{b}LU},
\end{split}
\end{equation}
where $Ro$ and $\delta$ represent the Rossby number and the Mach number, respectively. $\mu'$, $\nu'$ are the nondimensional viscosities. We assume $r_{p}/\rho_{h,m}=\delta$, which means that the ratio of the characteristic density perturbation $r_{p}$ to the mean $\rho_{h,m}$ is of order of $O(\delta)$. Thus, by noting that $\rho_{h,m}=\rho_{b}^{\gamma-1}$ we also have $r_{p}/\rho^{\gamma-1}_{b}=\delta$. Taking the new variables given in \eqref{eq9143} and \eqref{eq9142} into \eqref{CPE}, and dropping the primes we end up with the following dimensionless compressible primitive equations:
\begin{equation}\label{e3}
\begin{cases}
\partial_{t}\rho+\mathrm{div}_{h}\,(\rho u)+\partial_{z}(\rho w)=0,\\
\rho(\partial_{t}u+u\cdot\nabla_{h} u+w\cdot\partial_{z}u+\frac{1}{Ro}u^{\bot})+\frac{1}{\gamma\delta^{2}}\nabla_{h}\rho^{\gamma}=\mathbb{S}u,\\
\partial_{z}\rho^{\gamma-1}=-\theta,
\end{cases}
\end{equation}
where
\begin{equation*}
\rho=(1-\theta z+\delta r)^{\frac{1}{\gamma-1}},~~\partial_{z}r=0,~~\mathrm{and}~~\mathbb{S}u=\mu\triangle_{h}u+\mu\partial_{z}^{2}u+\nu\nabla_{h}\mathrm{div}_{h}\,u.
\end{equation*}
Using the new coordinates, the domain is correspondingly transformed to
\begin{equation*}
x=(x_{1},x_{2})\in\mathbb{T}^{2}:=[0,2\pi]\times[0,2\pi],~~0<z<\kappa:=\frac{L_{z}}{L},
\end{equation*}
where $\kappa$ is known as the aspect ratio parameter. Since $\rho$ is always non-negative, it's natural that
\begin{equation*}
\inf\limits_{z}\,(1-\theta z)=1-\theta\kappa\geq0,
\end{equation*}
thus $\theta\kappa\leq 1$.

\end{appendix}

\end{document}